\documentclass{article}

\usepackage{epsfig}
\usepackage{calc}
\usepackage{amssymb}
\usepackage{amstext}
\usepackage{amsmath}
\usepackage{color}

\usepackage{algorithm}
\usepackage{algorithmic}

\newcommand{\R}{\mathbb{R}}
\newcommand{\U}{\mathbb{U}}
\newcommand{\M}{\mathbb{M}}

\newcommand{\Z}{\mathbb{Z}}
\newcommand{\norm}[1]{\|#1\|}

\newcommand{\pb}[1]{$({\cal P}_{#1})$}
\newcommand{\pbb}[2]{$({\cal P}_{#1}(#2))$}
\newcommand{\abs}[1]{\left| #1\right|}

\newcommand{\Cnl}{{\cal C}_{NL}}
\newcommand{\Chyb}{{\cal C}_{\cal H}}
\newcommand{\Csubhyb}{\underline{\cal C}_{\cal H}}
\newcommand{\target}{{\cal T}}

\newtheorem{theorem}{Theorem}
\newtheorem{proposition}{Proposition}

\newtheorem{definition}{Definition}

\newtheorem{property}{Property}
\newtheorem{remark}{Remark}
\newtheorem{lemma}{Lemma}

\newenvironment{proof}{{\it Proof.}}{\hfill $\Box$\\}

\newenvironment{prob}[2]{\pb{#1}{\em #2}}{\\}

\begin{document}

\title{An hybrid system approach to nonlinear optimal control problems}
\author{Jean-Guillaume DUMAS\footnote{Universit\'e de Grenoble,
    laboratoire Jean Kuntzmann, Math\'ematiques Appliqu\'ees et Informatique, UMR CNRS 5224, Tour IRMA, BP 53, 38041 Grenoble, France. Email: Jean-Guillaume.Dumas@imag.fr} \hspace{1cm} Aude RONDEPIERRE\footnote{Universit\'e Paul Sabatier, Institut de Math\'ematiques de Toulouse, 118 route de Narbonne, 31062 Toulouse, France. Email: rondep@mip.ups-tlse.fr}}
\date{}
\maketitle

\begin{abstract}
We consider a nonlinear ordinary differential equation and want to
control its behavior so that it reaches a target by minimizing a cost
function. Our approach is to use hybrid systems to solve this problem:
the complex dynamic is replaced by piecewise affine approximations
which allow an analytical resolution. The sequence of affine models
then forms a sequence of states of a hybrid automaton. Given a
sequence of states, we introduce an hybrid approximation of the nonlinear
controllable domain and propose a new algorithm computing a controllable,
piecewise convex approximation. The same way the nonlinear optimal control problem is replaced by an hybrid piecewise affine one. Stating a hybrid maximum principle suitable to our hybrid model, we deduce the global structure of the hybrid optimal control steering the system to the target.
\end{abstract}



\section{Introduction}
Aerospace engineering, automatics and other industries provide a lot
of optimization problems, which can be described by optimal control
formulations: change of satellites orbits, flight planning, motion
coordination (see e.g. \cite{Fierro} or \cite{Pesch} for more applications in aerospace industry). In general those optimal control problems are fully nonlinear ; since the years 1950-1970, the theory of optimal control has been extensively developed and has provided us with powerful results like dynamic programming \cite{Bellman} or the maximum principle \cite{Pontryagin}. A large amount of theory has been developed and resolutions are mainly numerical.
 
We here consider a dynamical system whose state is described by the solution of the following ordinary differential equation (ODE):
\begin{equation}
\left\{
\begin{array}{l}
\dot{X}(t) = f(X(t),u(t))\\
X(0)=X_0
\end{array}
\right.\label{syst:nonlinear}
\end{equation}
We present a hybrid (discrete-continuous) algorithm controlling the
system (\ref{syst:nonlinear}) from an initial state $X_0$ at time
$t=0$ to a final state $X_f=0$ at an unspecified time $t_f$. To reach
this state, we  allow the admissible control functions $u $ to take
values in a convex and compact polyhedral set $\U_m$ of $\R^m$, in
order to minimize the following cost-function:
\begin{equation}
J(X,u) = \int_0^{t_f} l(X(t),u(t))dt\label{cout}
\end{equation}

There exists two main approaches for the resolution of nonlinear optimal control problems (see e.g. \cite{Trelat:book} for a review of different solving methods).

A first class of solving methods is based on the Bellman dynamic principle and the characterization of the value function $J(X,u)$ in terms of viscosity solutions of a nonlinear Hamilton-Jacobi-Bellman (HJB) equation \cite{Crandall:Lions:83,Barles:94,Bardi:Capuzzo,Capuzzo:Dolcetta,Bertsekas}. The related numerical methods, called direct methods, often uses time and/or space discretizations and are really efficient in small dimension. Unfortunately, their cost in memory is in general too expensive in high dimension. 

The second class of methods is based on the Pontryagin Maximum Principle \cite{Pontryagin,Gamkrelidze}. This principle provides a pseudo-Hamiltonian formulation of optimal control problems and some necessary optimality conditions. In this context, the foremost numerical methods are the indirect shooting methods. These methods are very efficient in any dimension but require an {\it a priori} knowledge of the optimal trajectory structure.

Today's challenge is to reach a compromise between direct and indirect
methods in order to enable the handling of high dimensional problems with a satisfactory precision and without any pre-simulation to locate optimal trajectories.

In this paper, we propose a different approach: the use of hybrid
computing \cite{DellaDora:2001:issac}. The idea is to approach nonlinear
systems like (\ref{syst:nonlinear}) by simplified models which we can
analytically study. Basically, an analytical approach must allow to
improve approximations as is done by \cite{Girard,Neurone:2003:HSCC}:
the level of details allows to reach a compromise between quantitative
quality of the approximation and the computational time. The latter
has been done e.g. for biological systems, where simplifications (in
relation to real data and in regard of model simulations) are
possible, see \cite{Neurone:2003:HSCC}. 
We here choose to use a hybrid
system modeling i.e. to approximate the original system
(\ref{syst:nonlinear}) by a continuous and piecewise affine one.
Contrary to \cite{Hedlund:Rantzer:99}, we want to focus on {\em both} the
discrete and continuous aspects of the hybrid systems. 
We thus first compute on the fly a polyhedral mesh of phase space {\em and} the control space and approach the global system by piecewise affine dynamics in
each mode: in each cell of this mesh, the system is piecewise affine
($\dot X(t)=AX(t)+Bu(t)+c$) and can be locally solved, mostly with
analytical tools.
We then find inside the cells the affine feedback controls \`a la
\cite{Habets:2006:PAHSS} which enable a global controllability. 
Finally, to optimize the choice of the solutions, 
we derive from \cite[Theorem 2.3]{Riedinger:Iung:Kratz:03} a hybrid
maximum principle where the transitions are not constrained. 
%
\nocite{Long:2006:IJRNC,Cho:2001:IJRNC}
\nocite{Baoticy:2003:EJC} 

The paper is organized as follows: \S\ref{sec:hybrid:approximation}
presents algorithms for the hybrid
approximation of nonlinear system. Next, in section
\ref{sec:controllability}, the controllability of the system is
studied and a new efficient algorithm computing a piecewise convex
approximation of the controllable set is developed. The last section
\S\ref{sec:optimality} is devoted to a heuristic approaching 
the resolution of optimal control
problems: we first validate our hybrid approach for this resolution 
and establish a suitable hybrid maximum
principle. We then show that the analysis of the optimality conditions
gives us the general structure of optimal trajectories.

\section{Hybrid Approximation of Nonlinear
  Systems}\label{sec:hybrid:approximation}
In this section, we model nonlinear control
systems by way of piecewise affine systems. Our main objective here is to replace the nonlinear system (\ref{syst:nonlinear}) by another one, complex enough to reproduce the intrinsic properties of the initial system, but also simple enough to be analytically studied.

Let us consider the general nonlinear control system
$\dot{X}(t) = f(X(t),u(t))$, 
where the admissible controls are measurable. Moreover
we assume that the bounded functions $u$ take values in a compact
polytope $\U_m$ of $\R^m$. 
Classical numerical methods like Euler or Runge-Kutta usually propose
time discretizations of nonlinear systems and then obtain
discretizations of solution trajectories. On the contrary, we here use
hybrid computation introduced in \cite{DellaDora:2001:issac} and whose
main idea is to compute an approximation of the nonlinear field $f$ by
means of a {\em state discretization}. Extended to control systems, the
principle is as follows: for a given mesh of the control and state
space $\R^n\times \U_m$, we compute a piecewise affine approximation
of the nonlinear vector field $f$. This partition, associated to the
so computed piecewise affine dynamic, defines a hybrid system which we
can thus analyze. We thus first describe the construction of a
piecewise affine approximation of the nonlinear system. We then
propose a complete study of the interpolation error and the
convergence of the approximation. Then we develop algorithms building
an implicit mesh of the state and control space $\R^n\times
\U_m$. This will enable us to define
our hybrid model for nonlinear control systems.

\subsection{Piecewise affine approximation of nonlinear systems}\label{hybridization}
In this paragraph we want to approximate the nonlinear system (\ref{syst:nonlinear}) by a system:
\begin{equation}
\dot{X}(t)=f_h(X(t),u(t))\label{syst:hybride}
\end{equation}
where $f_h$ is a piecewise affine approximation of the nonlinear vector field $f$.

\subsubsection{Computation of the piecewise affine approximation}
There are many possible linearizations by part~; especially in the
context of nonlinear systems without terms of control,
\cite{Chien:Kuh:1977} and \cite{Girard} propose a linearization
technique based on the interpolation of the nonlinear vector field at
the vertices of a given simplicial mesh of the state space. This
technique is extended in this paragraph to control dynamics like
(\ref{syst:nonlinear}) and therefore enables the computation of a
piecewise affine approximation with respect to the state $X$ {\em and} the control $u$.

Let $\Delta=(\Delta_i)_{i \in I}$ be a given simplicial mesh of $\R^n \times \U_m$. Any affine function in dimension $n+m$ is uniquely defined by its values at $(n+m+1)$ affine independent points. Thus, as done e.g. in \cite{Girard:2002:casc}, in each cell $\Delta_i$,  we can compute an affine approximation $f_h$ of $f$ by interpolation at the vertices of the cell. Moreover, only one supplementary function evaluation can be sufficient to compute an approximation when the system is evolving to an adjacent cell.

Then, we consider a simplicial cell $\Delta_i$ of $\Delta$, defined as
the convex hull of its $(n+m+1)$ vertices:
$\sigma_1,\dots,\sigma_{n+m+1}$. Let $f_i$ be the affine approximation
of $f$ computed by interpolation at the vertices of $\Delta_i$.
 We have that: 
$f_i(X,u)=A_i
X+B_i u+c_i$.
The approximation $f_i$ is computed by interpolation of $f$ at the
vertices of the cell $\Delta_i$, so that: $\forall j, \
f(\sigma_j)=f_i(\sigma_j)$. Hence, this proves that 
$f(\sigma_{j}) = [~A_i~|~B_i~]\sigma_j+c_i$ which in turn induces that
$\forall j=1,\dots,n+m+1,~f(\sigma_j) - f(\sigma_1) = [~A_i~|~B_i~](\sigma_j-\sigma_1)$.
We thus define $M_i$ to be the $(n+m) \times (n+m)$ matrix, whose
columns are the vectors of vertices: $\{\sigma_j-\sigma_1 ;
j=2,\dots,n+m+1\}$ and $F_i$ is the $(n+m) \times n$ matrix, whose
columns are the vectors
$\{f(\sigma_j) - f(\sigma_1) ; j=2,\dots,n+m+1\}$ so that $F_i=[~A_i~|~B_i~] M_i$.

Consequently, by linear independence of the vertices of the simplex $\Delta_i$,
the square matrix $M_i$ is non singular, so that we obtain the
piecewise affine approximation of figure \ref{fh_fi} by\footnote{
the $c_i$ can be solved as well with any column as 
$\forall j=1,\dots,n+m+1,\ f(\sigma_j) - [~A_i~|~B_i~]\sigma_j = c_i$.}:
\begin{figure}[htbp]
\begin{center}
\includegraphics[width=0.4\textwidth]{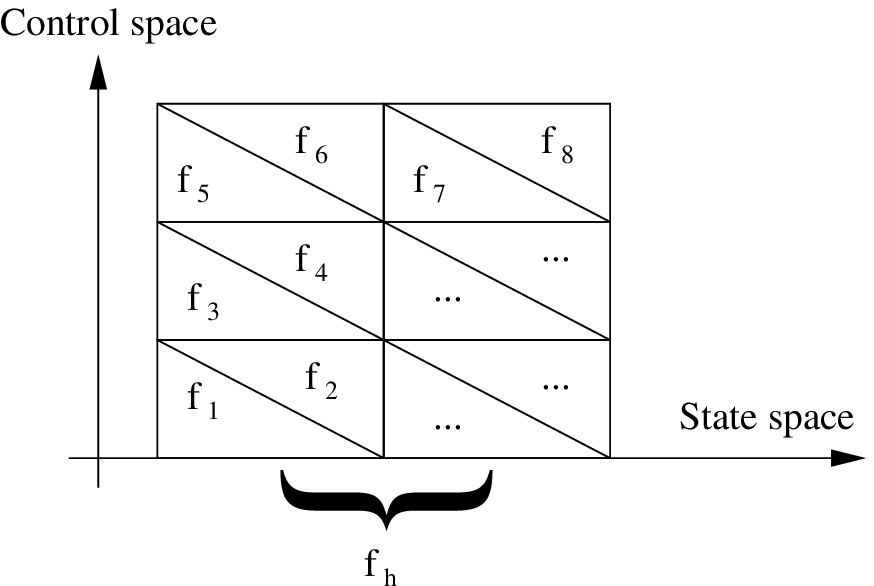}
\begin{equation}
\begin{array}{l}
\left\{\begin{array}{rcl}
[~A_i~|~B_i~] &=& F_i~M_i^{-1}\\
~\\
c_i &=& f(\sigma_1) - [~A_i~|~B_i~] \sigma_1\\
\end{array}\right.\\
~\\
f_h(X,u) =  A_i X + B_i u + c_i, ~~~\text{if}~ (X,u)\in \Delta_i
\end{array}
\label{approxhybride}
\end{equation}
\end{center}
\caption{Definition of the piecewise affine approximation}\label{fh_fi}
\end{figure} 

\subsubsection{Interpolation error and convergence of the Hybridization scheme}\label{sssec:convergence}
%
For autonomous nonlinear systems: $\dot{x}(t)=g(x(t))$, \cite[\S 11.1]{Girard} presents a study of the approximation error inside a compact set $D$ under  some regularity assumptions on the field $g$. In the context of nonlinear control systems $\dot{X}(t)=f(X(t),u(t))$, the control is in general only measurable so that the required regularity for $f$ is lost. 
Regularity assumptions on the initial nonlinear field $f$ are usually
chosen to fulfill the conditions of the Cauchy-Lipschitz theorem
applied to control systems, see
e.g. \cite{Sontag,Trelat:book,Rondepierre:phD}. From now on we therefore
assume that $f$ is continuous, Lipschitz continuous over any compact
set in $X$, uniformly in $u$. In other words, we have the condition $(A_1)$:
$$\begin{array}{cc}
(A_1) & \forall \Omega \Subset \R^n,~\exists
  L_{\Omega}>0,~\forall (X_1,X_2) \in \Omega^2,~\forall u \in \U_m,\\
& \norm{f(X_2,u)-f(X_1,u)}\leq L_{\Omega}\norm{X_2-X_1}.
\end{array}$$
Then for any given control $u$, we can prove the existence and uniqueness
of solution trajectories of the nonlinear system (\ref{syst:nonlinear}):
we first establish the regularity properties of the piecewise
affine approximation $f_h$, required by the Cauchy-Lipschitz
theorem. This then enables us to get an accurate bound on the error produced by our linearization scheme. We end this study by convergence results on our scheme.\\

Let us introduce the size $h$ of the simplicial mesh $\Delta=(\Delta_i)_{i\in I}$ defined by: 
$h=\sup_{i\in I} ~h_i $ where $h_i=\sup_{x,y\in \Delta_i} ~\norm{x-y}$
and $\norm{.}$ denotes the $\infty$-norm on $\R^{n+m}$. We first need
the continuity of the vector field $f_h$ which is an application
of \cite[proposition 11.1.3]{Girard} to a state-control space:
\begin{proposition}
$f_h$ is continuous over $\R^n\times \U_m$ and locally Lipschitz in $X$, i.e. for all compact subset $\Omega$ in $\R^n$, the restriction of $f_h$ to $\Omega\times \U_m$ is Lipschitz continuous in $X$, uniformly in $u$.\label{prop:property:fh}
\end{proposition}

%
Then, from this continuities of $f$ and $f_h$, we deduce the following bound on the magnitude of $f_h$:
\begin{lemma}
$f$ and $f_h$ are locally bounded in $X$, uniformly in $u$ on $\R^n\times \U_m$, i.e.:
$$\forall R>0,\exists C_R >0,\forall (X,u)\in B(0,R)\times \U_m,\norm{f_h(X,u)}\leq C_R$$
Moreover, if $f$ is bounded on $\R^n\times \U_m$, then $f_h$ is also bounded on $\R^n\times \U_m$ (with the same bound).\label{lem:bounded} 
\end{lemma}

%
Consequently for any given control $u$, the conditions of the
Cauchy-Lipschitz theorem are satisfied and we therefore obtain the
existence and uniqueness of solutions to the Cauchy problem: $\dot{X}(t)
= f_h(X(t),u(t)),~X(0)= X_0$. Now under some non restrictive
assumptions on the nonlinear vector field $f$, our interpolation error
can be evaluated with respect to the size $h$ of the mesh $\Delta$ of
$\R^n\times \U_m$:

\begin{proposition}[Interpolation error]
\noindent ~i.~If $f$ is L-Lipschitz on $\R^n\times \U_m$, then:
$$\sup\limits_{(X,u)\in \R^n \times \U_m} \norm{f(X,u)-f_h(X,u)} \leq \frac{4L(n+m)}{n+m+1}h= \varepsilon(h)$$
\noindent ii.~If $f$ is continuously differentiable on $\R^n \times \U_m$, then for any compact set $\Omega\subset \R^n$, then there exists $L_\Omega >0$ such that:
$$\sup\limits_{(X,u)\in \Omega\times\U_m}\norm{f(X,u)-f_h(X,u)} \leq \frac{4L_\Omega(n+m)}{n+m+1}h = \varepsilon_\Omega(h)$$
where: $L_{\Omega}=\max\limits_{q'\in I/\Delta_{q'}\cap(\Omega\times \U_m)\neq\emptyset}\left(\sup\limits_{(X,u)\in \Delta_{q'}} \norm{Df}\right)$.
\label{prop:interpolation:error}
\end{proposition}
The principle of the proof of proposition
\ref{prop:interpolation:error} is first to evaluate the interpolation
error inside one cell $\Delta_{q'}$ of the mesh $\Delta$. If the so
computed error does not depend on the given cell $\Delta_{q'}$, then
we deduce the global interpolation error ; otherwise the error is
evaluated on any compact set $\Omega\times \U_m$ of the state and
control space. See \cite[proposition 2.3.1]{Rondepierre:phD}
for more details.\\

Now, let $X_0\in \R^n$ be a given initial point and let the state $X(.)$
(resp. $X_h(.)$) be the corresponding solution of the initial
nonlinear (resp. piecewise affine) system: $\dot{X}=f(X,u)$
(resp. $\dot{X_h}=f_h(X_h,u)$). Now by introducing the attainable sets
${\cal A}(X_0,t)$ and ${\cal A}_h(X_0,t)$ i.e. the sets of points that
can be reached from $X_0$ by the related system in time $s\leq t$, we
obtain the following convergence results  \cite[proposition
2.3.2]{Rondepierre:phD} of our piecewise affine approximation:
\begin{proposition}[$C^1$ Convergence] ~\\
\noindent ~i.~If $f$ is Lipschitz continuous in $X$, uniformly in $u$, then for all $t>0$ such that $X(t)$ and $X_h(t)$ are well defined:
\begin{eqnarray}
\norm{X(t)-X_h(t)}&\leq& \frac{\varepsilon(h)}{L}(e^{Lt}-1)\nonumber\\
\norm{\dot{X}(t)-\dot{X}_h(t)}&\leq& \varepsilon(h)e^{Lt}.\nonumber
\end{eqnarray}
\noindent ii.~Let $f$ be Lipschitz continuous over any compact set $\Omega$ in $X$, uniformly in $u$ (hypothesis A1) and assume that the attainable sets ${\cal A}(X_0,s),~s>0$ are compact. Then for $t>0$,
\begin{eqnarray}
\norm{X(t)-X_h(t)}&\leq& \frac{\varepsilon_\Omega(h)}{L_\Omega}(e^{L_\Omega t}-1)\nonumber\\
\norm{\dot{X}(t)-\dot{X}_h(t)}&\leq& \varepsilon_\Omega(h)e^{L_\Omega t},\nonumber\end{eqnarray}
where $\Omega$ is defined by: $\Omega = {\cal A}(X_0,t)\cup{\cal A}_h(X_0,t)$.\label{prop:convergence}
\end{proposition}
We can thus see that the key point is the definition of the mesh $\Delta$ and its size $h$. We now show that meshing the whole space is not mandatory. We rather compute the mesh on the fly.

\subsection{Implicit simplicial mesh}\label{mesh}
Previously we have seen how to build a piecewise affine approximation
of the nonlinear dynamic $f$ for a given simplicial mesh of the state
and control domain. In order to perform our hybrid approximation, we need now a method to build a simplicial mesh of $\R^n\times \U_m$. As $\R^n\times \U_m$ is not bounded, it is algorithmically inconceivable to mesh the whole space $\R^n \times \U_m$. Our approach is then to implicitly define a mesh of our space so that the simplicial subdivision is made on the fly.\\

Let $h>0$ be the discretization step. There exists several ways to mesh a given domain. To each mesh corresponds only one piecewise affine approximation (\ref{syst:hybride}) by interpolation. Among all the possible meshes of the space $\R^n\times \U_m$, we choose to select those for which $(0,0)$ is a fixed point of the system (\ref{syst:hybride}). Indeed, if the target is $X_f=0$, we want to stay there once we have reached it. Therefore, $u(t_f)$ also has to be zero. A simple way to thus enforce $f_h(0,0)=(0,0)$ is to request the following property:
\begin{property}
If $0\in \Delta_i$ then $0$ is a vertex of $\Delta_i$.\label{property:0}
\end{property}
\begin{lemma}
Let us assume that $f(0,0)=0$. Let $\Delta$ be a mesh of $\R^n\times \U_m$ that satisfies the property \ref{property:0}. Then $(0,0)$ is also a fixed point of the piecewise affine system (\ref{syst:hybride}) built over $\Delta$.
\end{lemma}
\begin{proof} At the vertices of the mesh, $f(\sigma)=f_h(\sigma)$ by
  the interpolating constraints and $(0,0)$ is a fixed point of the
  initial system by definition. Therefore, $f_h(0,0) = f(0,0) = 0$.
\end{proof}

\begin{figure}[htbp]
\begin{center}
\begin{tabular}{c}
\includegraphics[width=0.3\textwidth]{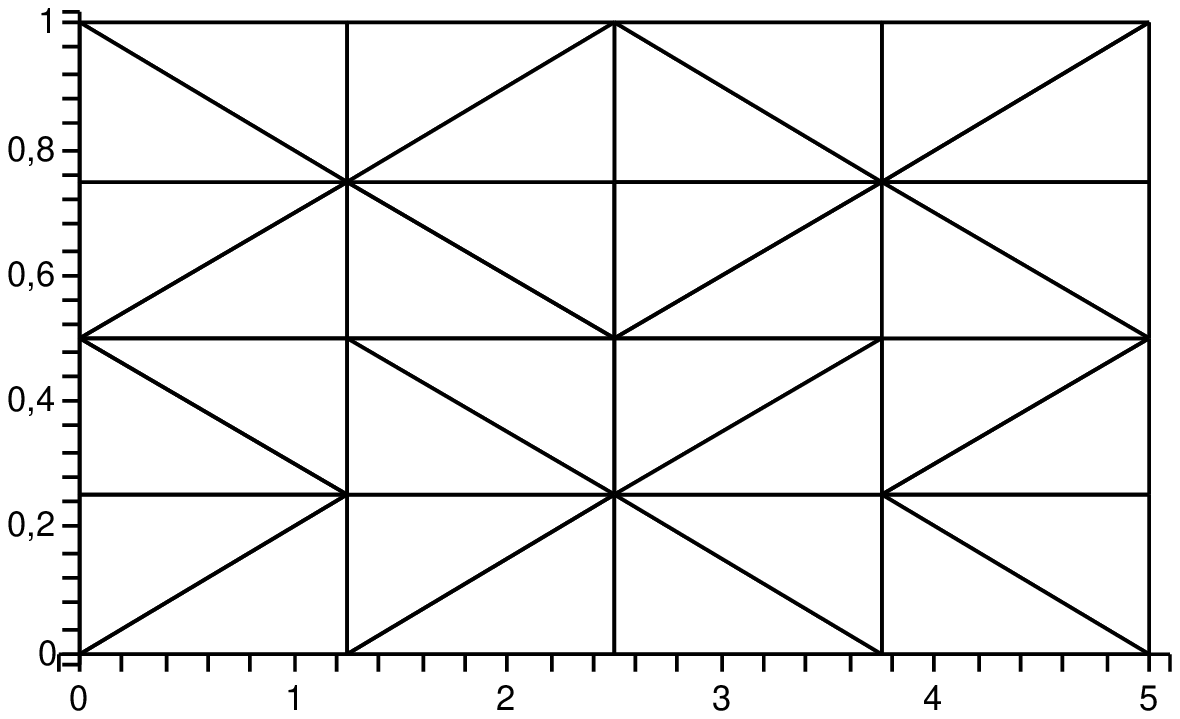}\\ 
(a) \\
\includegraphics[width=0.3\textwidth]{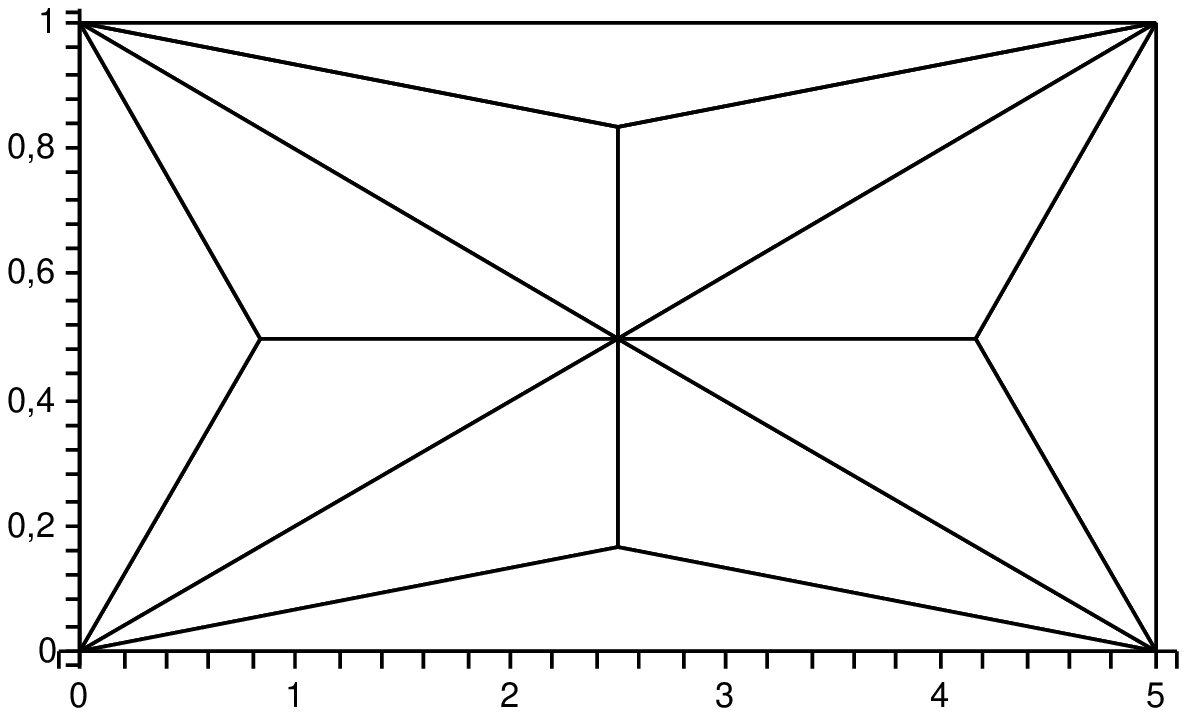}\\
(b)
\end{tabular}
\end{center}
\caption{Examples of meshes, for $n=m=1$, with (a) satisfying the
properties \ref{property:0} and \ref{property:mesh:2} and where (b) satisfies
property \ref{property:0} and not property \ref{property:mesh:2}.}\label{fig:mesh}
\end{figure}

Moreover we assume that the mesh $\Delta$ also satisfies the property:
\begin{property}
Let $D=(D_q)_{q\in {\cal Q}}$ be the projection of $\Delta$ over the state space $\R^n$.
\begin{enumerate}
\item $D$ is a simplicial mesh of $\R^n$.
\item Let $(i,j)\in I^2$. We state: $D_i=p_{\R^n}^{\bot}(\Delta_i)$
and $D_j=p_{\R^n}^{\bot}(\Delta_j)$. Then we have:
$D_i=D_j \mbox{ or } \mathring{D_i}\cap \mathring{D_j}=\emptyset$.
\end{enumerate}\label{property:mesh:2}
\end{property}
The property \ref{property:mesh:2} is a geometrical constraint on the
position of the cells $\Delta_i$ in $\R^n\times\U_m$. Indeed, we now
assume that, to each cell $D_i$ corresponds a ``column'' of cells
$\Delta_j$ whose projection over $\R^n$ is exactly $D_i$ (see for
example, figure \ref{fig:mesh}a). We note ${\cal K}(q) =  \{i\in I~;~p^{\bot}_{\R^n}(\Delta_{i})=D_q\}$ the set of all the indices of cells $\Delta_j$ belonging to the same column $q$.
Then each domain $D_q\times\U_m$ can be decomposed into exactly $card~{\cal K}(q)$ simplices in this way:
\begin{lemma}
$D_q\times \U_m=\bigcup\limits_{q'\in {\cal K}(q)}\Delta_{q'}$.\label{decomposition:mode:q}
\end{lemma}
The latter property will allow us in section \ref{sec:optimality} to have the same output constraints from a state cell $D_q$ for each $\Delta_i$ in the column defined by ${\cal K}(q)$.\\

We are now able to define an algorithm computing an implicit mesh of the space $\R^n\times \U_m$. To satisfy property \ref{property:mesh:2}, we start by computing an implicit mesh $D=(D_q)_{q\in{\cal Q}}$ of the state space. Next, we compute an explicit triangulation of the control domain $\U_m$ (small and bounded). We lastly deduce an implicit simplicial mesh $(\Delta_i)_{i\in I}$ of $\R^n \times \U_m$.

\subsubsection{Simplicial mesh of $\R^n$} 
The state space $\R^n$ is implicitly cut into n-dimensional cubes. Each cube is then meshed into $n!$ simplices and the resulting partition is a mesh of $\R^n$ (see \cite{Freudenthal},\cite[chapter 11]{Girard} for more details). Here, a $n$-cube $C$ is defined by one point $a=(a_i)_{i=1\dots n} \in \R^n$ and the length\footnote{The Euclidean norm is used to compute lengths or sizes in the mesh.} $h$ of its edges: $C = [a_1,a_1+h]\times \dots \times [a_n,a_n+h]$. This will also be denoted by $C= a+[0,h]^n $. In the same way we define a mesh of $\R^n$ into $n$-cubes:
\begin{definition}[\bf Mesh of $\R^n$ into $n$-cubes]
Let $a=(a_i)_{i=1\dots n}$ be a given point in $\R^n$ and $h>0$. Then $(a+kh+[0,h]^n)_{k=(k_1,\dots,k_n)\in \Z^n}$ is a mesh of $\R^n$ into $n$-cubes.\label{def:mesh}
\end{definition}
Now, let us consider a $n$-cube $C=[a_1,a_1+h]\times \dots \times [a_n,a_n+h]$. We then introduce ${\cal S}_n$ the set of permutations of $\{1,\dots,n\}$. For all $\varphi \in {\cal S}_n$, $D_{a,\varphi} = \{(x_1,\dots,x_n) \in \R^n ; 0 \leq x_{\varphi(1)} - a_{\varphi(1)}\leq \dots \leq  x_{\varphi(n)} - a_{\varphi(n)} \leq h\}$ is a simplex in $\R^n$, whose vertices are defined by: 
\begin{equation}
\left\{ \begin{array}{ll}
\forall i = 1,\dots,p, &x_{\varphi(i)}= a_{\varphi(i)}\\
\forall i = p+1,\dots,n, &x_{\varphi(i)}= a_{\varphi(i)} + h\\
\end{array}\right., p=0,\dots,n\label{def:vertex}
\end{equation}
Note that property \ref{property:0} is satisfied if $a=0$. Now each $n$-cube can be meshed into $n!$ simplices \cite{Freudenthal}:
\begin{proposition}
$(D_{b,\varphi})_{\varphi \in S_n}$ is a mesh of size $\sqrt{n}h$ of the $n$-cube $C = b +[0,h]^n$.\label{prop:Freudenthal}
\end{proposition}
We then check that the respective meshes of two adjacent $n$-cubes
coincide at their intersection so that we deduce a global mesh of
$\R^n$ as in \cite[proposition 3.1.2]{Rondepierre:phD}:
\begin{proposition}
$(D_{k,\varphi})_{k\in\Z^n,\varphi\in {\cal S}_n}$ is a simplicial mesh of size $\sqrt{n}h$ of the state space $\R^n$.\label{prop:state:mesh}
\end{proposition}

\subsubsection{Triangulation of $\U_m$} 
The control domain $\U_m$ is a bounded convex polyhedron, defined as
the convex hull of its vertices. The study of simplicial subdivisions
of such polytopes has been extensively developed in recent years and
provides us with some efficient tools to compute them (e.g. Delaunay
triangulation of the {\it Qhull} \cite{Qhull} software\footnote{tool
  for computation of convex hulls, Delaunay triangulation, Voronoi
  diagrams in 2d, 3d or higher dimension.}). To build a regular mesh of size $h$ of $\U_m$, we use an algorithm combining a Sierpinski-like discretization with the Delaunay triangulation.
\begin{remark}
The triangulation of the control domain is computed once and for all at the beginning of the algorithm: indeed, in every control problem, there is in general no initial condition on the control. Therefore in each state cell, we need to compute all the possible related control cells and the triangulation can not be computed on the fly as for the state. This is possible for the control polytope as it is bounded and usually much smaller than the state space. Also this same explicit triangulation will be reused in every column cell of our hybrid automaton.
\end{remark}

\subsubsection{Simplicial mesh of $\R^n \times \U_m$}
Let $(D_q)_{q\in {\cal Q}}$ and $(U_j)_{j\in I}$ be the respective
simplicial subdivisions of $\R^n$ and $\U_m$. The triangulation
$(\Delta_i)_i$ of the state and control space $\R^n \times \U_m$ is
built without new vertices via the Delaunay triangulation. This last
criterion guarantees the property \ref{property:mesh:2} to be
satisfied and gives the following size of the resulting mesh:
\begin{lemma}
Let $D_q$ and $U_j$ be two polytopes respectively in $\R^n$ and $\R^m$.
Then we have:
$$diam(D_q\times U_j) = \sqrt{diam(D_q)^2+diam(U_j)^2}$$
\end{lemma}
Combining this lemma with proposition \ref{prop:state:mesh} gives that our mesh
$\Delta$ is of size $\sqrt{n+1}h$.\\

In this section, we have described the main steps of the construction of an implicit simplicial mesh. Next, we will see that this mesh and the corresponding linear approximations, actually define a hybrid approximation of the initial system (\ref{syst:nonlinear}).

\subsection{Hybrid Automaton}\label{sec:hybrid:automaton}
%
\subsubsection{Definition of the hybrid model}\label{sec:automaton}
The state of a hybrid system is described by two variables: the first one, denoted by $q(t)$, is discrete and takes values in a countable set ${\cal Q}$. The second one, denoted by $X(t)$, is continuous and takes values in the space $\R^n$. In each discrete mode $q\in {\cal Q}$, the variable $X$ evolves continuously according to a dynamic $f_q$ specific to the mode $q$. Therefore the hybrid dynamic and consequently the solving methods, depend on the choice of the discrete set ${\cal Q}$.
Since the control functions are generally measurable and
not continuous (at best piecewise continuous), trajectories $(X(.),u(.))$ in the space $\R^n\times \U_m$ are a priori discontinuous whereas the trajectories $X(.)$ are continuous in the state space.
\begin{figure}[htbp]
\begin{center}
\scalebox{.5}{\input{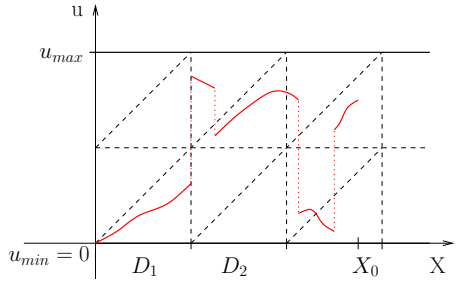}}
\caption{Example of discontinuous trajectory in the space $\R\times [0,u_{max}]$}\label{fig:trajectoire}
\end{center}
\end{figure}
Furthermore, as illustrated on figure \ref{fig:trajectoire}, the notion of transitions between two cells of the mesh $\Delta$ has no meanings with respect to the control. Indeed, an admissible control function can switch at any time to a given cell $\Delta_i$ to another $\Delta_j$, not necessarily adjacent to $\Delta_i$, and without reaching one boundary of the system.
In consequence the choice of the continuous variable, describing the
state of the hybrid system, is restricted to the state variable $X\in
\R^n$. We then introduce the following notations:
\begin{enumerate}
\item ${\cal Q}$ the countable set of indices of a simplicial mesh $(D_q)_{q\in {\cal Q}}$ of the state space $\R^n$.
\item ${\cal E}= \{(q,q') \in {\cal Q} \times {\cal Q}~;~\partial D_q \cap \partial D_{q'} \neq \emptyset\}$ the transition set: a transition between two discrete modes $q$ and $q'$ is allowed if the associated cells $D_q$ and $D_{q'}$ are adjacent.
\item ${\cal D}=\{D_q~;~q\in {\cal Q}\}$ the collection of simplicial cells of the mesh of $\R^n$. By construction, each domain $D_q$ induces a set of affine constraints on the state $X$ in the mode $q$:
$D_q=\{X\in \R^n~;~N_qX+L_q\leq O\}$.
\item ${\cal U} = \{\U_m~;~q \in {\cal Q}\} \subset \U_m$ the collection of control domains. Here, in each mode $q$ the control domain is unchanged.
\item ${\cal F} = \{f_q~;~q \in {\cal Q}\}$ the collection of piecewise affine field vectors defined by interpolation of the nonlinear field vector $f$ (see section \ref{hybridization}): 
$$\begin{array}{lccl}
\forall q\in {\cal Q},~f_{q}: &D_q\times\U_m &\rightarrow& \R^n\\
& (X,u) &\mapsto& f_h(X,u)
\end{array}
$$
According to lemma \ref{decomposition:mode:q}, the vector $f_q$ is locally defined by:
$\forall q'\in {\cal K}(q),\forall (X,u)\in \Delta_{q'}, f_q(X,u) = A_{q'}X+B_{q'}u+c_{q'}.$
\item ${\cal G} = \{G_e~;~e \in {\cal E} \}$ the collection of 
guards: 
$\forall e=(q,q') \in {\cal E} , G_e = \partial D_q \cap
\partial D_{q'}$.
\end{enumerate}
\begin{remark}
The definition of hybrid systems proposed at first by M.S. Branicky and al. in \cite{Branicky}, suggests the introduction of the collection ${\cal R} = \{R_e~;~e \in {\cal E} \}$ of Reset functions here defined by: 
$\forall e=(q,q') \in {\cal E}, \forall X \in G_e, R_e(X) = \{X\}$.
Indeed in our model, at any transition $e=(q,q')\in {\cal E}$, the corresponding cells $D_q$ and $D_{q'}$ are adjacent so that we do not need to reinitialize the continuous variable $X$.
\end{remark}
%
We can now define our hybrid model and the associated notion of
solution, extending that of \cite[\S 2.1]{Girard}:
\begin{proposition}
${\cal H= (Q,E,D,U,F,G)}$ is a hybrid system approximating the nonlinear system (\ref{syst:nonlinear}).
\end{proposition}
\begin{definition}
A solution (or execution) of the hybrid system ${\cal H}$ is a hybrid trajectory $\chi = (\tau,q,X)$ where:
\begin{enumerate}
\item[$\bullet$] $\tau=([t_i,t_{i+1}])_{i=0\dots r}$ ($t_i\leq
t_{i+1}$ if $i<r$) is a sequence of real intervals.
\item[$\bullet$] $q:\tau \rightarrow {\cal Q}$ describes the discrete state of ${\cal H}$
and $X:\tau \rightarrow \R^n$ the continuous behavior of ${\cal H}$.
\end{enumerate}
\noindent and satisfying, for all $i\in\{0,\dots,r\}$ such that $t_i<t_{i+1}$,
\begin{enumerate}
\item[$\bullet$] For $t \in ]t_i,t_{i+1}[$, $q(t)$ is constant
and $X(t)$ remains inside $D_{q(t)}$.
\item[$\bullet$] there is a measurable function $u:]t_i,t_{i+1}[\rightarrow \U_m$ such that:
$\forall t\in ]t_i,t_{i+1}[,~\dot{X}(t)=f_{q(t)}(X(t),u(t))$.
\end{enumerate}\label{def:hybrid:execution} 
\end{definition}
At any time $t$, the state of the hybrid system ${\cal H}$ is described by its position $X(t)$ in the state space $\R^n$ and the associated discrete mode $q(t)=q\in {\cal Q}$. For each nonempty time interval $]t_i,t_{i+1}[$, the system ${\cal H}$ is in mode $q(t)=q_i$ at the position $X(t)\in D_{q_i}$. In this mode, there exists an admissible control $u$ such that ${\cal H}$ evolves continuously according to the dynamic: $\dot{X}(t)=f_{q_i}(X(t),u(t))$. As soon as the system reaches a guard $G_{(q_i,q_{i+1})}$ at time $t_{i+1}$, a discrete transition $e=(q_i,q_{i+1})\in{\cal E}$ may occur from the mode $q_i$ towards the mode $q_{i+1}$ and we start again the same process from the position $X(t_{i+1})$ in mode $q_{i+1}$ at time $t_{i+1}$. Note that $t_i$ and $t_{i+1}$ are the respective entering and exiting times in the cell $D_{q_i}$.
%
%
%
\subsubsection{Computation of hybrid data on the fly}\label{ssec:local:data}
Let $X_0$ be the current position in the state space $\R^n$ of the
hybrid system ${\cal H}$ defined in the previous section. In this
section we first present an explicit algorithm to compute on the fly the current
discrete mode associated to the position $X_0$ i.e. the cell $D_q$
which contains the point $X_0$. We will assume in
particular that the hybrid automaton is not Zeno (no infinite
number of switch can occur in a finite time, see e.g. \cite[\S
2.3]{Girard} and \cite{Zeno:HS}).
Then we want to compute the local control constraints induced by each cell $\Delta_{q'}$ in the mode $q$.
\paragraph{Computation of the current state}
The principle of the algorithm is based on the definition of the implicit mesh $\Delta$ of the state and control space $\R^n\times\U_m$ (see section \ref{mesh}): we first compute the $n$-cube and then the simplex inside this cube that contain the given position $X_0$.

\noindent{\it Computation of the $n$-cube $C$ such that $X_0\in C$:}   according to the definition \ref{def:mesh}, the $n$-cube $C$ is defined by one origin point and the length $h$ of its edges. Here the origin of the mesh $\Delta$ has to be chosen at the origin $0$ of $\R^n\times\U_m$, so that the searched cube $C$ can be written as: $C = kh + [0,h]^n$, where $kh=(k_1h,\dots,k_nh)\in h\Z^n$. Since we want to compute all the cubes $C$ to which $X_0=(x_1,\dots,x_n)\in \R^n$ belongs, we necessarily obtain:
$$\forall i\in\{1,\dots,n\}, \left\{\begin{array}{ll}
k_i = \left[\displaystyle\frac{x_i}{h}\right] &\mbox{ if } x_i \notin h\Z\\
\\
k_i \in \{\displaystyle\frac{x_i}{h}-1,\frac{x_i}{h}\} &\mbox{ otherwise}
\end{array}\right.$$
\noindent{\it Computation of the simplex $D_q$ such that $X_0\in D_q$}: let $C=kh +[0,h]^n$ be a $n$-cube containing $X_0$. According to the proposition \ref{prop:Freudenthal}, $C$ can be divided into $n!$ simplices defined by:
\begin{eqnarray}
\forall \varphi\in {\cal S}_n,D_{\varphi}=\{X\in \R^n~;~0\leq X_{\varphi(1)}-k_{\varphi(1)}h \leq \dots\label{def:Dphi} \\
\dots \leq X_{\varphi(n)}-k_{\varphi(n)}h\leq h\}\nonumber
\end{eqnarray}
To avoid the explicit computation of all the $D_\varphi$ contained in
the cube $C$, the trick is to directly compute the permutation
$\varphi$ which defines the searched simplex:
\begin{enumerate}
\item Form $L:=[x_1-k_1h,\dots,x_n-k_nh]$;
\item Apply a sort by selection algorithm to the list $L$. We so obtain the sorted list $\widetilde L$ and the permutation $\varphi$ that performs the sort on the list $L$.
\item Deduce the simplex $D_\varphi$ satisfying: $X_0\in D_\varphi$
\end{enumerate}

\begin{figure}[htbp]
\begin{center}
\includegraphics[width=0.2\textwidth]{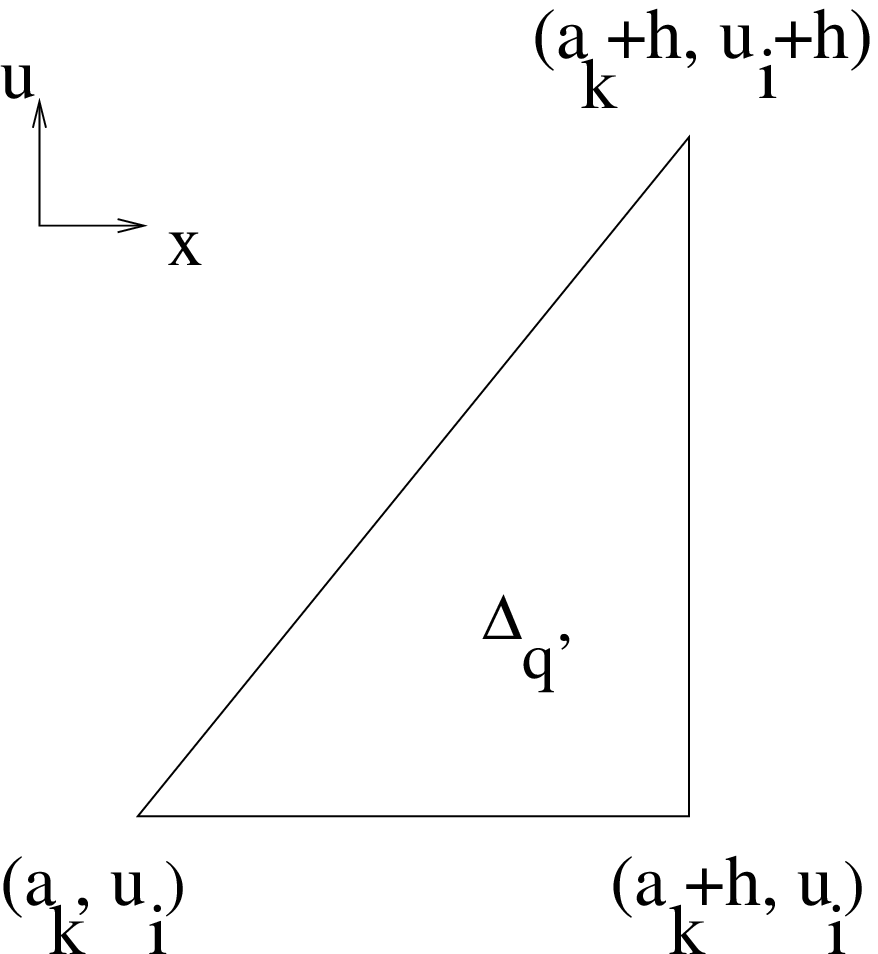}

$(x,u) \in \Delta_{q'}$\\ 
$\Updownarrow$\\ 
$x \in D_q = [a_k,a_k+h]$\\ 
and \\
$u \in U_{q'}(x) = \{u \in \U_1 ; 0 \leq u-u_i \leq x-a_k\}$
\end{center}
\caption{Definition of $D_q$ and $U_{q'}(x)$ in cell $\Delta_{q'}$ for
$n=m=1$}\label{fig:U}
\end{figure}

\paragraph{Computation of the local constraints on the control}
We consider a given cell $D_q$ in the state space and again 
${\cal  K}(q)=\{q'\in I~/~p^{\bot}_{\R^n}(\Delta_{q'})=D_q\}$
is the set of cells in $\Delta$ whose projection in $\R^n$ is
exactly $D_q$.
By construction, in each cell $\Delta_q$, the control depends on the
state. In this paragraph we want to compute the explicit control
constraints. For that matter, we define $ U_{q'}(X)$ 
to be the set of control constraints in the cell
$\Delta_{q'}$ for a given $X\in \R^n$:
$\forall X\in D_q,~\forall q'\in {\cal K}(q),~U_{q'}(X) = \{u \in \R^m ~/~(X,u)\in \Delta_{q'}\}$.
The problem now is to determine the geometrical structure of
$U_{q'}(X)$ in $\R^m$.
For instance, on figure \ref{fig:U}, we can remark that the
control domain $U_{q'}(x)$ is a segment (namely an 1-simplex) in $\R$
when $x\in ]a_k,a_k+h]$, and is reduced to one point $\{u_i\}$, when
$x=a_k$.

Let us take $X\in D_q$ and $q' \in {\cal K}(q)$. By definition,
$\Delta_{q'}$ is a simplex in $\R^{n+m}$, so that it can be defined by
a system of $(n+m+1)$ affine independent inequalities: $MY+d \geq 0$ where $M\in \M_{n+m+1,n+m}(\R)$ and $d\in \R^{n+m+1}$. We now define the left and right parts of $M$ as follows:
$$M_1=M\left[\begin{array}{c}
I_n\\
0
\end{array}\right]\hfill M_2=M\left[\begin{array}{c}
0\\
I_m
\end{array}\right]$$
so that: $M=[~M_1~|~M_2~]$. Thus we have $\forall (X,u)\in \Delta_{q'}, M_1X + M_2 u +d \geq 0$.
And we hence obtain a characterization of $U_{q'}$:
\begin{equation}\label{eq:uq}
\begin{array}{ll}
\forall X\in D_q,U_{q'}(X) = \{u \in \R^m /(X,u)\in \Delta_{q'}\}\\
\phantom{\forall X\in D_q,U_{q'}(X) }= \{u \in \R^m /M_2u+(M_1X+d)\geq 0\}
\end{array}
\end{equation}
Therefore $U_{q'}(X)$ is determined by a finite number of affine inequalities and is a polyhedral set. Moreover, since $\Delta_q$ is bounded by construction, it follows that $U_{q'}(X)$ is a bounded polyhedral set, i.e. a polytope. Furthermore, it also has the following property:
\begin{proposition}
For all $X\in \mathring{D_q}$. Then, for all $q'\in {\cal K}(q)$,
$U_{q'}(X)$ is a $m$-simplex in $\R^m$.\label{prop:Uq}
\end{proposition}
The proof is given in appendix \ref{apdx:uq} ; the idea is to count the number of intersections between the whole control space and the $n$-faces of the simplex $\Delta_{q'}$.

%

In this section, we have proposed a model of hybrid automaton for the
approximation of nonlinear control systems $\dot{X}=f(X,u)$. The
construction of the hybrid model is performed on the fly, so that at
any time we are able to determine the state and the associated local
control constraints. 
Next section addresses the controllability problem, namely the existence of (non compulsorily optimal) solutions via this hybrid approximation.

\section{Approximation of the controllable domain}\label{sec:controllability}
In this section, we focus on the controllability of the nonlinear system (\ref{syst:nonlinear}) towards a given polyhedral target. For a given initial point $X_0\in \R^n$ in the state space, we want to study the existence of admissible trajectories of the system that steers $X_0$ to the considered target at some finite time.
The theory of linear control systems without constraints on the control have been extensively developed and provides some powerful results like the Kalman controllability criterion \cite{Kalman:1963,Kalman:Falb:Arbib} for autonomous linear systems like $\dot{X}=AX+Bu$. The use of this criterion can be extended to nonlinear control systems by using local controllability techniques an so provides local controllability results \cite{Lee:Markus,Trelat:book}. Some global results can be obtained for a class of systems which depends affinely on the control \cite{Jurdjevic,Lobry:1970}.

Our goal is to propose a constructive approach of
nonlinear controllability by the study of the controllable domain
(the set of initial points that can be steered to the target). 
We thus present now methods and algorithms to approximate the
controllable domain of nonlinear systems by the way of the hybrid
model defined in section \ref{sec:hybrid:automaton}.
We first present our hybrid approach for the controllability of
nonlinear systems: the nonlinear controllable domain is approximated
by the hybrid one. Then, after stating some useful properties of the
hybrid controllable domain, we develop a new algorithm computing a
convex approximation of the hybrid controllable domain. Some
experimental results will then validate our algorithm.

\subsection{Hybrid approach to the nonlinear controllable domain}\label{par:hybrid:approach}
In section \ref{sec:hybrid:approximation}, we have built a hybrid
model approximating every nonlinear control system like
(\ref{syst:nonlinear}). We now need to verify that the use of our hybrid model enables the approximated solving of the initial control problem ; namely we analyze the approximation of the nonlinear controllable domain by the hybrid one.
Let us specify the notion of controllability and controllable domain to a given target $\target_f$:
\begin{definition}
$X_0 \in \R^n$ is controllable to $\target_f$ if and only if there
exists $T \geq 0$ and $u:[0,T]\rightarrow \U_m$ measurable, such that
the following system (S) admits a solution $X(.)$ over $[0,T]$:
$$ (S) \left\{\begin{array}{lcl}
\dot{X}(t) &=& f(X(t),u(t))\\
X(0) &=& X_0,~~X(T) \in \target_f
\end{array}\right.$$
The set of controllable points in $\R^n$ to $\target_f$ is called the {\em controllable domain} of the system (\ref{syst:nonlinear}).\label{def:controllability}
\end{definition}

This definition can be applied to any control system, and
particularly to the hybrid automaton previously described. 
Therefore, any $X_0\in \R^n$ is controllable
to the target $\target_f$ for the hybrid system ${\cal H}$ if and only
if there exists an admissible hybrid trajectory, in the sense of
definition \ref{def:hybrid:execution}, that steers $X_0$ to
$\target_f$ in some finite time $T>0$. In particular, for any
controllable point $X_0\in \R^n$, the hybrid system may accept several
executions that steer $X_0$ to the target, and several different cell
paths $q=(q_i)_i$ between $X_0$ and $\target_f$.
%

\subsubsection{Global error of the hybrid approximation}
Let ${\cal C}_{NL}(\R^n)$ and ${\cal C}_{\cal H}(\R^n)$ be the respective controllable domains of the nonlinear system (\ref{syst:nonlinear}) and its hybrid model ${\cal H}$ towards a given target $\target_f$. To validate the hybrid approximation, we need to evaluate the distance between ${\cal C}_{NL}(\R^n)$ and ${\cal C}_{\cal H}(\R^n)$, namely their Hausdorff distance defined by:
$$\begin{array}{l}
d_H(\Cnl(\R^n),\Chyb(\R^n)) \\
\phantom{blabla}= \max\left(\sup\limits_{X\in \Chyb(\R^n)}\inf\limits_{Y\in \Cnl(\R^n)} \norm{X-Y},\right.\\
\phantom{blablabblablabblablablala}\left.\sup\limits_{Y\in \Cnl(\R^n)}\inf\limits_{X\in \Chyb(\R^n)} \norm{X-Y}\right)
\end{array}$$

Let $X_0\in \Cnl(\R^n)$ be a controllable point to $\target_f$ by the nonlinear system (\ref{syst:nonlinear}). Hence, according to definition \ref{def:controllability}, there exist a finite time $T>0$ and an admissible control $u_0:[0,T]\rightarrow \U_m$ such that the system:
$\dot{X}(t) = f(X(t),u(t)),~~X(0)=X_0$
has a solution $X$ satisfying: $X(T)\in \target_f$. Then we can easily build a new point $Y_{h,0}$ that is controllable by the hybrid system ${\cal H}$ in the same time $T$ and control $u_0$ as for $X_0$:
$Y_{h,0} = \widetilde{Y}(T)$
where $\widetilde{Y}$ is the unique solution of the system: $\dot{\widetilde Y}(t) = -f_h(\widetilde Y(t),u_0(T-t)),~~\widetilde Y(0)=X(T)$ obtained from (\ref{syst:nonlinear}) by time reversal on the time interval $[0,T]$. We thereupon deduce: 
$Y_{h,0}\in \Chyb(\R^n) ~~\mbox{ and }~ \inf_{Y\in \Chyb(\R^n)} \norm{X_0-Y} \leq \norm{X_0-Y_{h,0}}$.
According to the convergence results established in proposition
\ref{prop:convergence}, and
assuming that  $f$ is Lipschitz in $X$, uniformly in $u$, we thus
obtain: $\norm{X_0-Y_{h,0}}\leq
\displaystyle\frac{\varepsilon(h)}{L}(e^{LT}-1)$, where
$\varepsilon(h)$ denotes the interpolation error between $f$ and its
piecewise affine approximation $f_h$ ;  we consequently have shown that for all $ X_0\in \Cnl(\R^n)$:
$$\exists T>0,\inf\limits_{Y\in \Chyb(\R^n)} \norm{X_0-Y} \leq \frac{\varepsilon(h)}{L}(e^{LT}-1).$$
Moreover, starting from a controllable point of the automaton
${\cal H}$, we prove in the same manner that:
$\forall Y_{h,0} \in \Chyb(\R^n),~\exists T_h\geq 0,~ \inf_{X\in \Cnl(\R^n)} \norm{X-Y_0} \leq \displaystyle\frac{\varepsilon(h)}{L}(e^{LT_h}-1).$
Since the final $T$ is a priori unbounded, we introduce the notion of
controllable domain in time less than or equal to $T>0$, denoted by
${\cal C}_{NL,[0,T]}(\R^n)$ in the nonlinear case and ${\cal C}_{{\cal
    H},[0,T]}(\R^n)$ in the hybrid one. We have shown:
\begin{proposition}
Let the nonlinear function $f:\R^n\times\U_m\rightarrow\R^n$ be Lipschitz continuous in $X\in \R^n$, uniformly in $u\in\U_m$. Then:
$$d_H\displaystyle\left({\cal C}_{NL,[0,T]}(\R^n),{\cal C}_{{\cal H},[0,T]}(\R^n)\displaystyle\right) \leq \varepsilon(h)\displaystyle\frac{e^{LT}-1}{L},$$
where the interpolation error $\varepsilon(h)$ satisfies: $\lim\limits_{h\rightarrow 0} \varepsilon(h)=0$.\label{prop:CD:approximation}
\end{proposition}
Using the second part of the proposition \ref{prop:convergence}, it is
then straightforward to generalize this result as:
\begin{proposition}
Let $f:\R^n\times\U_m\rightarrow \R^n$ be Lipschitz continuous over any compact set in $X\in \R^n$, uniformly in $u\in\U_m$. Then there exists $L_{\Omega}>0$ such that:
$$d_H\displaystyle\left({\cal C}_{NL,[0,T]}(\R^n),{\cal C}_{{\cal H},[0,T]}(\R^n)\displaystyle\right) \leq \varepsilon_{\Omega}(h)\displaystyle\frac{e^{L_{\Omega}T}-1}{L_{\Omega}},$$
where $\Omega\subset\R^n$ is the compact set defined by the union of ${\cal C}_{NL,[0,T]}(\R^n)$ and ${\cal C}_{{\cal H},[0,T]}(\R^n)$.\label{prop:CD:approximation:2}
\end{proposition}

\subsubsection{Exploration reduction: controllable domain in finite paths of modes}\label{ssec:exploration:reduction}
We have shown that the hybrid controllable domain is a valid approximation of the nonlinear one. From now on, we consider that a given point $X_0$ in the state space is approximatively controllable by the initial nonlinear system (\ref{syst:nonlinear}) to the considered target if it belongs to the controllable domain of its hybrid approximation.
Unfortunately the controllable domain of a hybrid system, like its
attainable set, are rarely computable. One of the main reasons is that
the number of discrete transitions needed to reach some modes of
the hybrid system is unbounded \cite[chapter 7]{Girard}: by time reversal, the
computation of the controllable domain comes down to the computation
of the attainable sets from the target. Then, if the given point $X_0$
is ``far from the target'', all the admissible paths of cells between
$X_0$ and the target have to be explored. Hence an infinite number of
discrete transitions are to be simulated and the computation can not
be done in finite time. Therefore an approach that is more favorable
to the simulation is required. We first assume that:
$$\begin{array}{cl}
(A_2) & \mbox{\it The target $\target_f$ is a subset of some cell $D_q$ in the}\\
& \mbox{\it state space $\R^n$: $\exists q\in {\cal Q},\target_f\subset D_q$}
\end{array}$$
The assumption $(A_2)$ is non restrictive: indeed, if there are several cells $D_q$ whose intersection with the target $\target_f$ is non-empty, then $\target_f$ can be decomposed into the union of its intersections with each cell:
$\target_f =\bigcup_{q\in {\cal Q}/D_q\cap\target_f\neq
  \emptyset} \left(D_q\cap\target_f\right).$
There, we have several targets to reach (as many new targets as
intersections) each of those satisfying the hypothesis $(A_2)$.

We first consider a given mode $q\in {\cal Q}$ of our hybrid model ${\cal H}$ and a local target $\target_q$ (i.e.: $\target_q\subset D_q$). We introduce the concept of controllable domain in the given mode $q$:
\begin{definition}[Controllable domain in a mode $q$]
Let $q\in {\cal Q}$ be a given discrete mode of our hybrid model ${\cal H}$ and $\target_q$ a subset of $D_q$. $X_0\in D_q$ is controllable to $\target_q$ in mode $q$ if there exist a finite time $T>0$ and an admissible control $u$ such that the hybrid problem:
$$\dot{X}(t)=f_h(X(t),u(t)),~~ X(0)=X_0,~~ X(T)\in \target_q$$
has a solution, denoted by $X_h$, satisfying: $\forall t\in [0,T],~X_h(t)\in D_q$. 

The controllable domain in mode $q$ is the set of all the controllable points in mode $q$ and is denoted by $\Chyb(\target_q,D_q)$.\label{def:DC:modeq}
\end{definition}
Hence, any point in $D_q$ from which the local target $\target_q$ is attainable while remaining within the cell $D_q$, necessarily reaches $\target_{q'}\supset\target_q$ ; hence:
\begin{lemma}
$\target_q \subset \target_{q'}\subset D_q \Rightarrow \Chyb(\target_q,D_q) \subset \Chyb(\target_{q'},D_q)$
\end{lemma}

\begin{figure}[htbp]
\begin{center}
\scalebox{.5}{\input{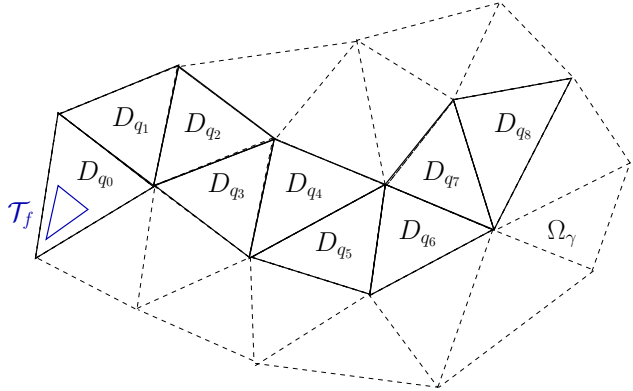}}
\caption{Definition of the domain $\Omega_\gamma$ along the path $\gamma=(q_i)_{i=0\dots 8}$}\label{fig:Omega:gamma}
\end{center}
\end{figure}

The definition \ref{def:DC:modeq} can be easily extended to any subset of the state space, and especially to the set $\Omega_\gamma$ associated to $\gamma=(q_i)_{i=1\dots r}$. This set, shown on figure \ref{fig:Omega:gamma}, is defined by 
$\Omega_\gamma = \bigcup_{i=1}^r D_{q_i}$
where $\gamma=(q_i)_{i=0\dots r}$ is a finite sequence of adjacent modes of our hybrid automaton ${\cal H}$ satisfying 
$\forall i=1,\dots,r-1,~D_{q_i}\cap D_{q_{i+1}} \neq \emptyset$.

\begin{figure}[htbp]
\begin{center}
\scalebox{.625}{\input{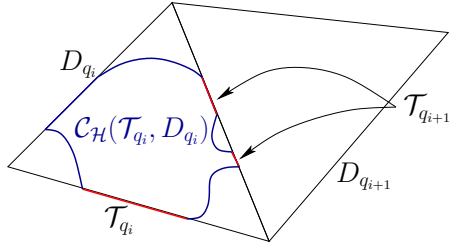}}
\caption{Spreading rule for the controllable domain from mode $q_i$ to mode $q_{i+1}$}\label{fig:propagation}
\end{center}
\end{figure}

Let $\target_f$ be a given subset of $\Omega_\gamma$ and the target
for the controllability problem. We introduce the set
$\Chyb(\target_f,\Omega_\gamma)$ of points that can be steered in
finite time by the hybrid model ${\cal H}$ to the target $\target_f$
while remaining within $\Omega_\gamma$, in the given sequence of modes.
Our idea is to build an approximation of the controllable set inside $\Omega_\gamma$, recursively defined along the path. Assuming without loss of generality that: ${\cal T}_f\subset D_{q_0}$, the principle consists in starting from the cell $D_{q_0}$: we first compute the local controllable domain in mode $q_0$. Then its intersection with the current guard $G_{(q_0,q_1)}$ defines a new target now in mode $q_1$ and we pursue the approximation the same way (the spreading rule is illustrated on figure \ref{fig:propagation}). The algorithm stops when the last cell $q_r$ is reached or when the intersection between the current controllable domain and the next guard is empty.

\noindent This induces a recursive definition of an under-approximation of the hybrid controllable domain inside $\Omega_\gamma$:
\begin{proposition}
We define:
$\Csubhyb(\Omega_\gamma)=\Chyb(\target_{q_0},D_{q_0}) \cup \Chyb(\target_{q_1},D_{q_1})\cup \dots \cup \Chyb(\target_{q_r},D_{q_r}),$
where $\Chyb(\target_{q_i},D_{q_i})$ denotes the set of points inside $D_{q_i}$ from which the local target $\target_{q_i}\subset D_{q_i}$ is attainable in finite time while remaining within the cell $D_{q_i}$. The current target $\target_{q_i}$ is recursively defined by:
$$\left\{\begin{array}{ll}
\tau_0=\target_f&\\
\tau_{i+1}=\Chyb(\tau_i,D_{q_i})\cap G_{(q_i,q_{i+1})},&
i=0,\dots,r-1 \end{array}\right.$$
$\Csubhyb(\Omega_\gamma)$ is an under-approximation of $\Chyb(\Omega_\gamma)$, hence~: $\Csubhyb(\Omega_\gamma)\subset \Chyb(\Omega_\gamma)$.\label{prop:under:approximation}
\end{proposition}
\begin{proof} We proceed by induction.
Let us define for $i=0,\dots,r$, the property ${\cal
  P}(i)~:~\Chyb(\target_{q_0},D_{q_0}) \cup \dots\cup
\Chyb(\target_{q_i},D_{q_i})\subset \Chyb(\Omega_\gamma)$.

First, let $X_0\in \Chyb(\target_{q_0},D_{q_0})$. Since $D_{q_0}\subset \Omega_\gamma$ and $\target_{q_0}=\target_f$, the point $X_0$ is actually controllable to the $\target_f$ while remaining within $\Omega_\gamma$. This proves that $X_0\in \Chyb(\Omega_\gamma)$ and ${\cal P}(0)$ is true.

Second, let us assume that the property ${\cal P}(i)$ is satisfied and $X_0\in \bigcup_{k=0}^{i+1} \Chyb(\target_{q_k},D_{q_k})$. By assumption, if $X_0\in \bigcup_{k=0}^i\Chyb(\target_{q_k},D_{q_k})$ then $X_0\in \Chyb(\Omega_\gamma)$. We can thus assume $X_0\in \Chyb(\target_{q_{i+1}},D_{q_{i+1}})$.
By definition of the local controllable domain inside the cell
$D_{q_{i+1}}$, there exists $Y_0\in \target_{q_{i+1}}$ such that the
hybrid system ${\cal H}$ steers $X_0$ to $Y_0$ while remaining inside
$D_{q_{i+1}}$, and henceforth inside $\Omega_\gamma$. Moreover by construction:
$\target_{q_{i+1}} = G_{q_i,q_{i+1}}\cap \Chyb(\target_{q_i},D_{q_i}),$
so that $Y_0$ is steerable to the target $\target_f$ while remaining inside $\Omega_{\gamma}$. By transitivity, we conclude that $X_0$ is controllable to the target $\target_f$ while remaining within $\Omega_\gamma$: $X_0\in \Chyb(\target_f,\Omega_\gamma)$ and the property ${\cal P}(i+1)$ is satisfied.
\end{proof}

The decomposition of the controllable domain in $\Omega_\gamma$ in local controllable sub-domains expressed in proposition \ref{prop:under:approximation} draws the structure of the algorithm presented in section \ref{sec:convex:approximation} for the hybrid controllable domain computation. Actually, as shown next, it is also possible to use this spreading in different directions.

\subsubsection{Expansion of the controllable domain in $r$ iterations}\label{ssec:spreading}
\begin{figure}[htbp]
\begin{center}
\scalebox{.45}{\input{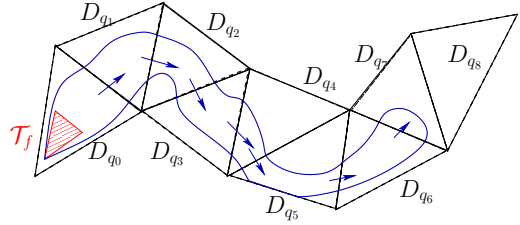}}\\
\caption{Controllable under-approximation in a given finite cells path $\gamma=(q_i)_{i=0,\dots,8}$.}\label{fig:underapprox1}
\end{center}
\end{figure}

\begin{figure}[htbp]
\begin{center}
\scalebox{.425}{\input{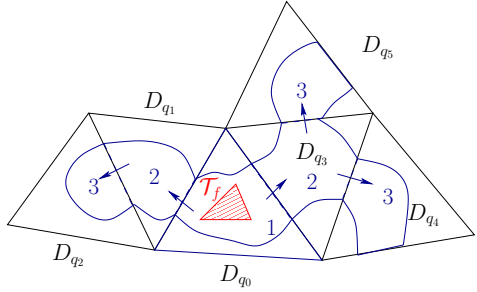}}
\caption{Controllable under-approximation in $3$ iterations from the target $\target_f$.}\label{fig:underapprox2}
\end{center}
\end{figure}

In paragraph \ref{ssec:exploration:reduction} we have proposed a
constructive approach of the hybrid controllable domain to a given target $\Omega_\gamma$. Its key points are described by the proposition \ref{prop:under:approximation}. In practise, this process may be used in two different manner illustrated on figures \ref{fig:underapprox1} and \ref{fig:underapprox2}.

The first is an extension of the process described previously: we are
given a finite sequence $\gamma=(q_i)_{i=0,\dots,r-1}$ of adjacent
modes of the hybrid model ${\cal H}$ and we compute the controllable
under-approximation as expressed in proposition
\ref{prop:under:approximation} and illustrated on figure
\ref{fig:underapprox1}. A weakness of this method is the arbitrary choice of the path in which we want to compute our controllable under-approximation.
As we will see later a way to exhibit such a path is to start with a
cell large enough to contain both the source and the target. Then such
an approximation at a given level could guide the selection of an
interesting path
at a more refined approximation level.

The second idea is based on a more global approach of the hybrid controllable domain. There, we compute an under-approximation of the global hybrid controllable domain $\Chyb(\R^n)$ by time reversal starting from the cell containing the target $\target_f$. The difference is that we only impose a maximum distance to the target. 
The principle is illustrated on figure \ref{fig:underapprox2}: we compute the cell $D_{q_0}$ of the state space containing the target ($\target_f\subset D_{q_0}$). Starting from $D_{q_0}$ by time reversal, we compute the local controllable domain to $\target_f$ in this cell and then its intersections with {\em all the guards (i.e. facets) of $D_{q_0}$}. Each non empty intersection defines a new target and a new cell at the next iteration. This process is pursued at the $p$-th iteration from the controllable domains computed at the $(p-1)$-th iteration.
%

\subsection{Computation of a convex approximation of the hybrid controllable domain}\label{sec:convex:approximation}
We have studied the approximation of the controllable
domain to a target $\target_f$ of a nonlinear system by the
controllable domain of its hybrid approximation. In view of
simulations, we also developed a more constructive approach for the
hybrid controllable domain computation. This approach is based on
local controllable domains propagation and reduces the state space
exploration.
Now, we want to compute the set of controllable points (or attainable
points by time reversal) in $\R^n$, i.e. the set of initial points
that can be steered to a given target $\target_f$ in finite time. 
The considered class of systems is the one of piecewise affine hybrid systems defined in section \ref{sec:hybrid:automaton}.

A first method for the attainable sets computation consists in
extracting, by abstraction, a simplified model equivalent to the
initial one for the regarded property
\cite{Alur:Dang:Ivanvic:02,Asarin:Dang:04,Tabuada:Pappas:01}. The main
problem is to find the right abstraction that enables to conclude
under a reasonable complexity. Another class of methods developed for
safety verification is based on approximations 
by unions of polytopes. For instance, \cite{Asarin:2003:HSCC}
compute an over-approximation of the attainable set and can thus certify that the system can not escape from an admissible set of states. 
On the contrary our idea is to guaranty as much as possible the
controllability of a given initial point. In \cite{issac:2005}, we
proposed an algorithm to compute an under-approximation in time $T>0$
of the controllable set without state constraints. In
\cite{icinco:2005}, an extension of this algorithm to piecewise affine
systems under constant control constraints is presented. Here, these techniques are adapted to our hybrid model ${\cal H}$ for bounded convex polyhedral targets. From now on, we assume:
$$\begin{array}{cl}
(A_3) & \mbox{\it The target $\target_f$ is a bounded polytope of $\R^n$ defined}\\
& \mbox{\it as the convex hull of its vertices.}
\end{array}$$

\subsubsection{Convex approximation within a given state cell}\label{sssec:cell:approximation}
Let us consider a given discrete mode $q\in {\cal Q}$ of our hybrid model, or equivalently a given cell $D_q$ in the state space. Let $\target_q\subset D_q$ be a local target in mode $q$ satisfying the hypothesis $(A_3)$:
$\target_q=Conv(s_{q,1},\dots,s_{q,k}).$
In this paragraph, we address the problem of approximating the local controllable set $\Chyb(\target_q,D_q)$ i.e. the set of points in $D_q$ steerable in finite time to $\target_q$ while remaining within $D_q$. Our main idea is to compute some specific controllable points in $D_q$ ; their convex hull defines a convex approximation of the hybrid controllable domain in mode $q$.

\begin{figure}[htbp]
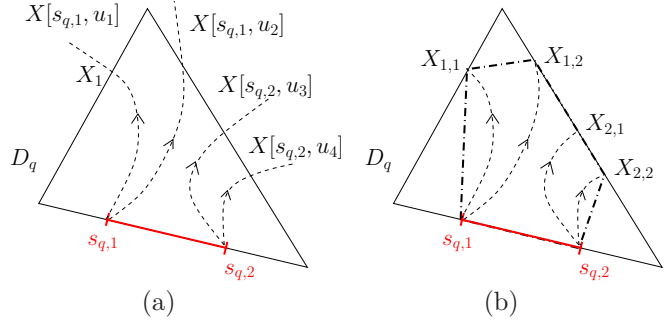

\begin{center}
\begin{tabular}{cc}
\scalebox{.425}{\input{Approxprinciple.pstex_t}} &\hspace{.2cm}
\scalebox{.425}{\input{ConvexApprox.pstex_t}}\\
(a) & (b)
\end{tabular}
\caption{Principle of the convex approximation in a state cell $D_q$: (a) along the trajectory $X[s_{q,1},u_1]$, only points between $s_{q,1}$ and $X_1$ are controllable to $s_{q,1}\in\target_q$ (b) Convex approximation}\label{fig:convexapprox:principle}
\end{center}
\end{figure}

Now, any point $X_0$ in $D_q$ is controllable to a target $\target_q$ in mode $q$ if it is attainable from $\target_q$ by time reversal. Consequently, a simple way to compute controllable points in a mode $q$ is: for an admissible control $u_j$, we compute by time reversal the hybrid trajectories $X[s_{q,i},u_j]$ coming from the vertex $s_{q,i}$ of $\target_q$. We then introduce their intersections, when they exist, with the boundary of $D_q$:
\begin{equation}
X_{i,j}=X[s_{q,i},u_j](T_i) \label{controllable:points}
\end{equation}
where $T_i = \sup\{t>0 ; X[s_{q,i},u](t) \in \partial D_q\}$. As illustrated on figure \ref{fig:convexapprox:principle}, when the trajectory $X[s_{q,i},u_j]$ remains inside $D_q$ on a non empty time interval, all points between $s_{q_i}$ and the corresponding intersection $X_{i,j}$ is controllable to the target $\target_q$. Our approximation will be the convex hull of all the points computed this way.\\

The main difficulty in this approach is the choice of control functions $u$. Indeed they have to be simple enough to enable the efficient computation of the related trajectories, but also significant enough to obtain a good quality approximation of the hybrid controllable domain.
According to the construction of the mesh $\Delta=(\Delta_{q'})_{q'\in I}$ of the state and control space developed in section \ref{mesh}, there exists a column of cells $\{\Delta_{q'};q'\in {\cal K}(q)\}$ in $\R^n\times \U_m$ whose projection on the state space $\R^n$ is exactly $D_q$. We recall that the cells of this column are defined by their indices this way:
${\cal K}(q) = \{i\in I~;~p^{\bot}_{\R^n}(\Delta_{i})=D_q\}$.
By definition of our hybrid model, in each cell $\Delta_{q'}$ of this column, the hybrid dynamic $f_h$ is affine. Our idea is then to choose a control function such that the control and state trajectory $(X(t),u(t))$ remains within only one cell $\Delta_{q'}$ ($q'\in {\cal K}(q)$) while $X(t)$ remains within the state cell $D_q$:
\begin{equation}\label{ss:cd}
\exists T>0,~\forall t\in [0,T],
X(t) \in D_q~\mbox{ and }~(X(t),u(t)) \in \Delta_{q'}
\end{equation}
We thus ensure that the hybrid dynamic is affine over the time
interval $[0,T]$. Our method is then to consider {\em feedback affine controls}:
\begin{equation}
u_j(t)=F_jX(t)+g_j\label{feedback:affine:control}
\end{equation}
such that $(X(t),u(t))$ evolves on the edges of our state and control mesh $\Delta$. Such controls at the position $X(t)$ correspond to the vertices of the control domains $U_{q'}(X(t))$ defined in paragraph \ref{ssec:local:data} and their definition is illustrated on figure \ref{fig:choix:controle}.\\
\begin{figure}[htbp]
\begin{center}
\scalebox{0.65}{\input{choixcontrole.pstex_t}}
\caption{Control $u_i(t)$ chosen for the computation of a convex approximation of the hybrid controllable domain in mode $q$}\label{fig:choix:controle}
\end{center}
\end{figure}

The principle of our approximation algorithm \ref{algo:convex:approximation} in mode $q$ is as follows: for each previously defined control $u_j$ and each vertex $s_{q,i}$ of the target $\target_q$,
\begin{enumerate}
\item Test if the trajectory by time reversal $X[s_{q,i},u_j]$ remains inside $D_q$ on a non empty time interval (notice that for this test, we do not need to explicitly compute the trajectory, see \cite[chapter 5, paragraph 5.1.2]{Rondepierre:phD}). If not, stop.
\item Computation of $X[s_{q,i},u_j]$ by solving:
$\dot{X}(t) = A_{q'}X(t)+B_{q'}u_j(t)+c_{q'}$, i.e. $\dot{X}(t) = (A_{q'}+B_{q'}F_j)X(t)+(B_{q'}g_j+c_{q'}).$
\item Computation of its first intersection with the boundary of $D_q$.
\end{enumerate}
\begin{algorithm}[htbp]
\caption{\texttt{\bf ConvexApproximation}: Convex approximation of the
  hybrid controllable domain in a mode}\label{algo:convex:approximation}
\begin{algorithmic}[1]
\REQUIRE ${\cal H}$, $q$ a discrete mode of the hybrid automaton ${\cal H}$, $\target_q = Conv(s_{q,1},\dots,s_{q,k})$ the local target in mode $q$.
\ENSURE A convex approximation of the controllable domain $\Chyb(\target_q,D_q)$.
\STATE {Initialization:} $\Lambda:=\emptyset$;\\
\STATE\COMMENT{Computation of the indices $q'\in {\cal K}(q)$ of cells $\Delta_{q'}$ whose projection is exactly $D_q$~:}\\
${\cal K}(q):=\{q' \in I~/~p^{\bot}_{\R^n}(\Delta_{q'})=D_q\}$
\STATE\COMMENT{Set of already tested controls:} $\Sigma:=\emptyset$.
\FORALL{$q' \in {\cal K}(q)$}
	\STATE\COMMENT{Computation of the polytope $U_{q'}(X)$ of control constraints (see \S\ref{ssec:local:data}):}\\
$U_{q'}(X):=ControlConstraints(\Delta_{q'})$;
\STATE\COMMENT{For each non tested vertex of $U_{q'}(X)$,}
   \FORALL{i from 1 to card($U_{q'}(X)$) such that $U_{q'}(X)[i] \notin \Sigma$}
       	\STATE\COMMENT{$i$-th vertex of $\Sigma$:} $u:=t\rightarrow
U_{q'}(X(t))[i]$; $\Sigma:=\Sigma \cup \{u\}$;
        \STATE\COMMENT{For each vertex $s_{q,j}$ of the target $\target_q$,}
        \FORALL{j from 1 to k}
        	\IF{the trajectory $X[s_{q,j},u]$ remains in $D_q$ on a non empty time interval}
		\STATE\COMMENT{Computation of the intersection of $X[s_{q,j},u]$ with $\partial D_q$:}\\
		$(X_f,ExitFace):=OutCell^\dagger({\cal H},q,s_{q,j},u)$;
		\STATE\COMMENT{$X_f$, if it exists, is controllable,}\\
		\STATE {\bf if} $X_f\neq \emptyset$ {\bf then} $\Lambda:=\Lambda \cup \{X_f\}$; {\bf end if}
	        \ENDIF
	\ENDFOR
  \ENDFOR
\ENDFOR
\STATE Return Conv($\Lambda$).
\end{algorithmic}
\hrule
\vspace{.1cm}
{\footnotesize $\dagger$ $(X_f,ExitFace):=OutCell({\cal H},q,X_0,u)$ computes the intersection point $X_f$ of the trajectory $X[X_0,u]$ of the hybrid automaton ${\cal H}$ with the state cell $D_q$. $ExitFace$ is the facet of $D_q$ satisfying: $X_f\in ExitFace$.}
\end{algorithm}

Let $\Lambda(\target_q,D_q)$ be the convex hull of the so computed controllable points $X_{i,j}$ defined by (\ref{controllable:points}) and computed by the algorithm \ref{algo:convex:approximation}. By construction, we deduce:
\begin{proposition}
$\Lambda(\target_q,D_q)$ is a convex approximation of the hybrid
controllable domain $\Chyb(\target_q,D_q)$ in mode $q$. Moreover if $\Chyb(\target_q,D_q)$ is convex, $\Lambda(\target_q,D_q)$ is a {\em guaranteed under-approximation}:
$$\Chyb(\target_q,D_q) \mbox{ convex }\Rightarrow \Lambda(\target_q,D_q)\subset \Chyb(\target_q,D_q).$$
\end{proposition}

\subsubsection{Expansion of the convex approximation and approximation error}
Let us consider a finite sequence of adjacent discrete modes of the hybrid automaton ${\cal H}$ satisfying $\forall i\in\{1,\dots,r-1\},~D_{q_i}\cap D_{q_{i+1}} \neq \emptyset $ and $\target_f\subset D_{q_0}$. As in section \ref{ssec:exploration:reduction}, we introduce the set $\Omega_\gamma$ associated to the path $\gamma$ and defined by: $\Omega_\gamma = \bigcup\limits_{i=0}^rD_{q_i}.$
According to proposition \ref{prop:under:approximation}, the hybrid controllable domain in $\Omega_\gamma$ can be under-approached by the union of local controllable domains successively defined in the cells $D_{q_i}$ of the path $\gamma$:
\begin{eqnarray}
\Csubhyb(\Omega_\gamma)=\Chyb(\target_0,D_{q_0}) \cup \Chyb(\target_1,D_{q_1})\cup \dots\phantom{blabla}\nonumber\\
 \cup \Chyb(\target_r,D_{q_r})\subset \Chyb (\Omega_\gamma)\label{controllable:domain:gamma}
\end{eqnarray}
where $\target_{q_i}$ denotes the local target in mode
$q_i$. Similarly, we define an approximation $\Lambda(\Omega_\gamma)$ of $\Csubhyb(\Omega_\gamma)$ by successive local convex approximations:
\begin{equation}
\Lambda(\Omega_\gamma)=\Lambda(\tau_0,D_{q_0}) \cup \Lambda(\tau_1,D_{q_1})\cup \dots \cup \Lambda(\tau_r,D_{q_r})\label{controllable:approximation:gamma}
\end{equation}
where $\Lambda(\tau_i,D_{q_i})$ denotes the convex approximation
computed by the algorithm \ref{algo:convex:approximation} 
and $(\tau_i)_{i=0\dots r}$ is the sequence of local targets defined along the path $\gamma$ by recurrence as follows:
\begin{equation}
\left\{\begin{array}{l}
\tau_0=\target_f\\
\tau_{i+1} = \Lambda(\tau_i,D_{q_i}) \cap G_{(q_i,q_{i+1})},~i=0,\dots , r
\end{array}\right.\label{cible:approchee}
\end{equation}
We easily check that the targets $\tau_i$, $i=0,\dots,r$ are polytopes (i.e. bounded and convex polyhedra).
Furthermore, each local target $\tau_i$ is a convex approximation of $\target_{q_i}$, so that $\Lambda(\tau_i,D_{q_i})$ is actually a convex approximation of the controllable domain $\Chyb(\target_{q_i},D_{q_i})$:
\begin{proposition}
For all $i=0,\dots,r$, the local targets $\tau_i$ defined by the
recursive relation (\ref{cible:approchee}), are respective convex
approximations of targets $\target_{q_i}$ and
$\Lambda(\tau_i,D_{q_i})$ is a convex approximation of order
$\sqrt{n}h$ of the controllable domain $\Chyb(\tau_i,D_{q_i})$ in mode $q_i$.\\
Moreover, if the sets $\Chyb(\target_{q_i},D_{q_i})$, $i=0,\dots,r$, are convex then:
$\forall i=0,\dots,r,~\tau_i\subset \target_{q_i}~~\mbox{ and } \Lambda(\tau_i,D_{q_i}) \subset \Chyb(\target_{q_i},D_{q_i}).$
\label{prop:global:approximation}
\end{proposition}

Now, we want to evaluate the approximation error, namely the Hausdorff distance $d_H(\Chyb(\Omega_\gamma),\Lambda(\Omega_\gamma))$ between the controllable domain $\Chyb(\Omega_\gamma)$ and its piecewise convex approximation $\Lambda(\Omega_\gamma)$. Since $\Csubhyb(\Omega_\gamma)$ is an under-approximation of $\Chyb(\Omega_\gamma)$ and according to the decompositions (\ref{controllable:domain:gamma}) and (\ref{controllable:approximation:gamma}), we deduce:
$$\begin{array}{lcl}
d_H(\Chyb(\Omega_\gamma),\Lambda(\Omega_\gamma)) & \leq & 
d_H\displaystyle\left(\Csubhyb(\Omega_\gamma),\Lambda(\Omega_\gamma)\right)\\
\\
&\leq & \max\limits_{i=0\dots r} d_H(\Chyb(\target_{q_i},D_{q_i}),\Lambda(\tau_i,D_{q_i}))
\end{array}$$
Moreover, for any cell $D_{q_i}$, $i=0,\dots,r$, the distance of the
local controllable domain $\Chyb(\target_{q_i},D_{q_i})$ to its convex
approximation $\Lambda(\tau_i,D_{q_i})$ is bounded by the size of the
mesh ${\cal D}$ of the state space. Hence, according to
\ref{prop:state:mesh}, we have: 
$\forall i=1,\dots,r,~d_H(\Chyb(\target_{q_i},D_{q_i}),\Lambda(\tau_i,D_{q_i}))\leq \sqrt{n}h$
and we can conclude by: $$d_H(\Chyb(\Omega_\gamma),\Lambda(\Omega_\gamma)) \leq \sqrt{n}h.$$
\subsection{Experience on the non linear spring}\label{ssec:exp}
\begin{figure}[htbp]
\begin{center}
\includegraphics[width=.25\textwidth]{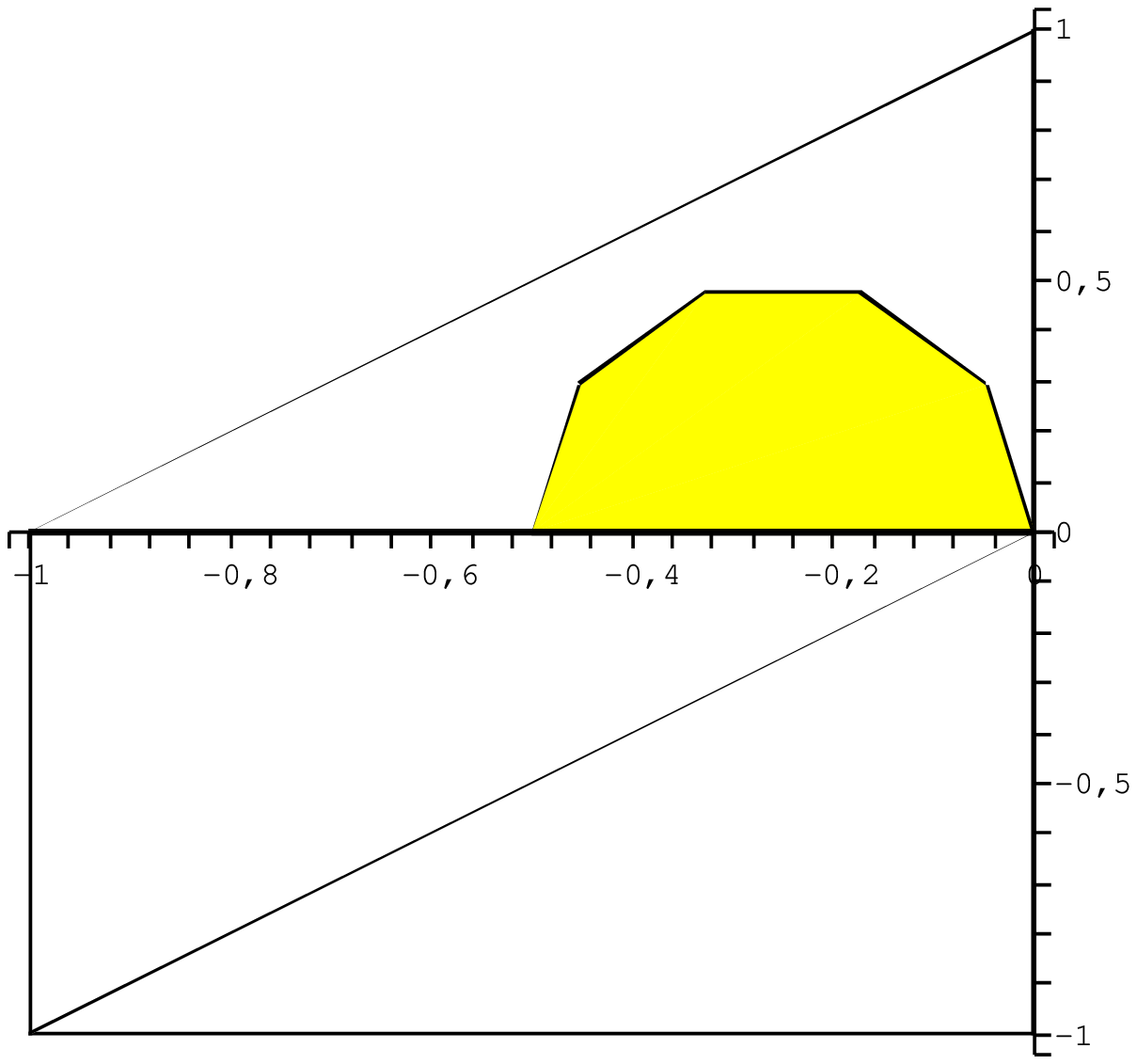}
\includegraphics[width=.25\textwidth]{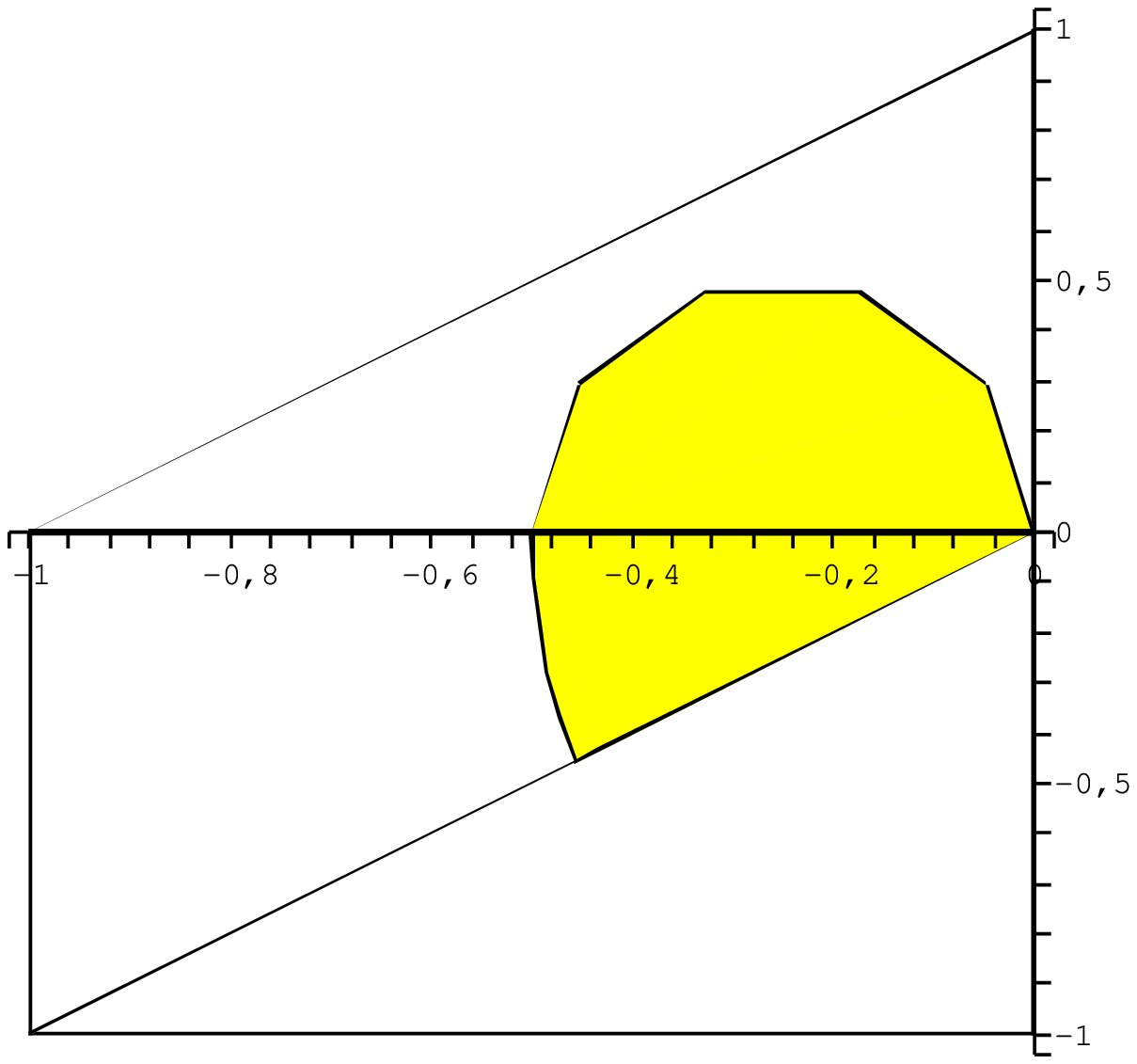}
\includegraphics[width=.25\textwidth]{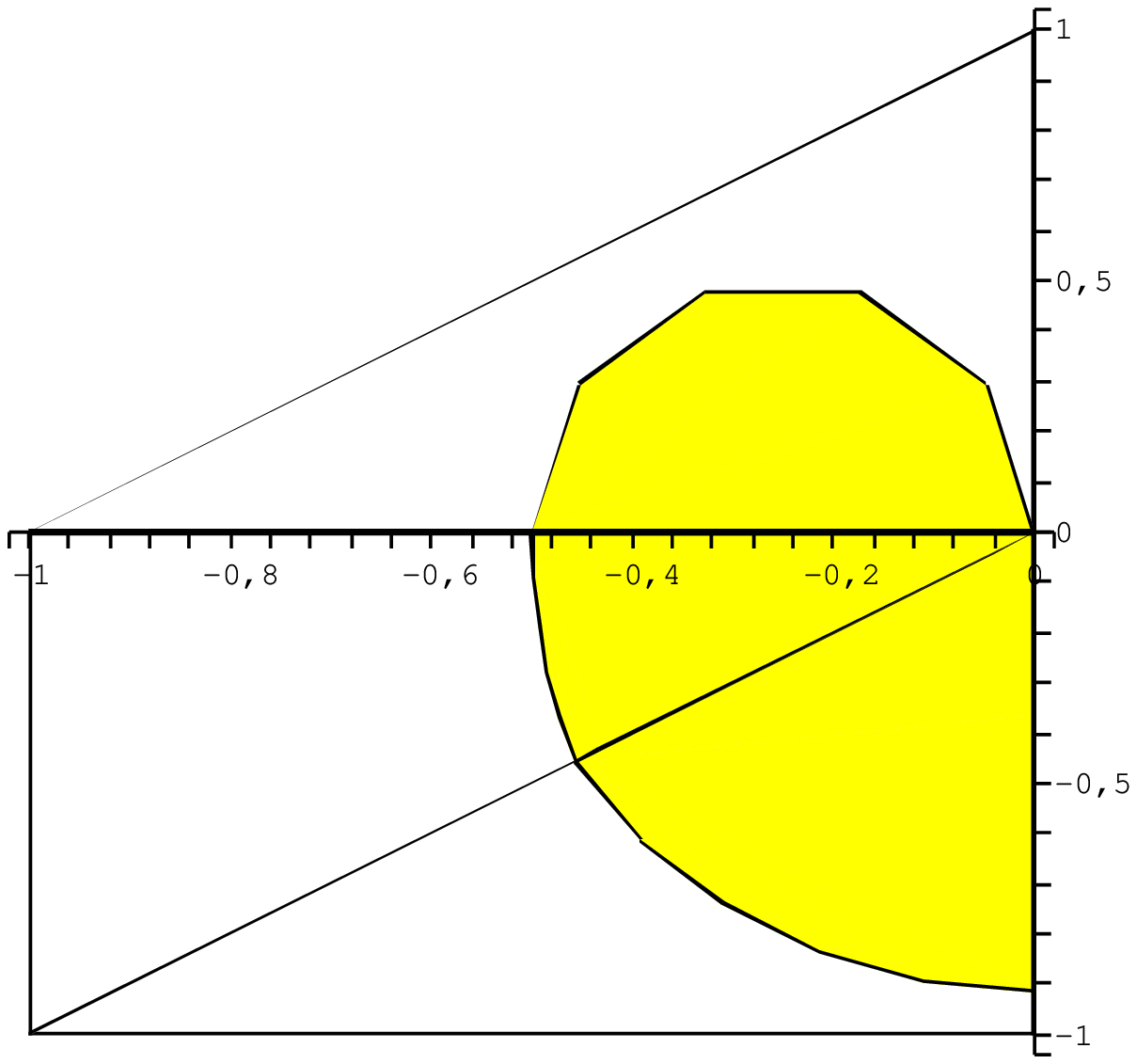}
\caption{Successive controllable approximations inside three adjacent cells for $h=1$.}\label{fig:ressort}
\end{center}
\end{figure}

\begin{figure}[htbp]
\begin{center}
\scalebox{.5}{\includegraphics{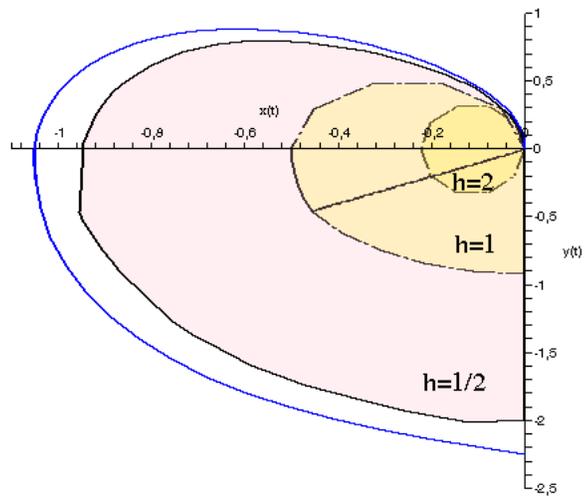}}
\caption{Convex controllable approximations at scales $h=2$, $h=1$ and $h=1/2$. Our approximations tend to the exact controllable domain delimited by the blue curve.}
\end{center}\label{fig:ressort:2}
\end{figure}

In this section, we consider the small example of the nonlinear spring presented in \cite[paragraph 7.3.1]{Trelat:book}:
\begin{equation}
\left\{\begin{array}{lcl}
\dot{x}(t)&=&y(t)\\
\dot{y}(t)&=&-x(t)-2x(t)^3+u(t)
\end{array}\right.
\end{equation}
where the control constraints are: $\forall t\geq 0,~|u(t)|\leq 1$. In this section, we want to apply the previously described algorithms to compute a controllable approximation of the set of initial positions $(x_0,y_0=\dot{x}_0$), from which the spring can be steered to its equilibrium position $(0,0)$.\\
In this example, we perform the convex approximation algorithm on three successive adjacent cells. Figure \ref{fig:ressort} shows the three successive hybrid controllable approximations for the discretization step $h=1$. Figure \ref{fig:ressort:2} shows the result of our algorithm at three different discretization level $h\in \{2,1,\frac{1}{2}\}$. The outside blue curve delimits the real nonlinear controllable domain over the three considered cells: as we can observe, the hybrid approximations tend to the exact nonlinear domain when $h$ decreases.


\subsection{Orbital transfer}
We now consider a satellite in orbit around the Earth. Our control problem is the minimum time transfer of this satellite from an initial orbit to a geostationary one. 

Let $m$ be the mass of the satellite and $u$ its thrust. The satellite dynamic is described by the equations:
\begin{eqnarray}
\ddot{r}=-\mu\frac{r}{\norm{r}^3}+\frac{u}{m},\label{satellite:mouvement}\\
\dot{m}=-\frac{\norm{u}}{v_e},~ v_e\in \R^+.\label{satellite:masse}
\end{eqnarray}
where $\mu$ denotes the gravitation constant. The control $u$ is assumed bounded: $\norm{u}\leq F_{max}$. The numerical values of the involved physical constants are given in Table \ref{table:ctes:phys}.\\

\begin{table}[htbp]
\begin{center}
\begin{tabular}{|ll|}
\hline
{\bf Parameter} & {\bf Value}\\
\hline
$\mu$ & $5165,8620912$ Mm$^3$.h$^{-2}$\\
$m_0$ & $1500$ kg\\
$F_{max}$ & 3 N\\
\hline
\end{tabular}
\caption{Physical constants for the orbital transfer problem}\label{table:ctes:phys}
\end{center}
\end{table}
The mathematical model of the satellite is a variable mass three dimensional model. Simplified models have been proposed such as coplanar transfer \cite{Caillau:Noailles:01} and/or constant mass models \cite{Bonnard:Caillau:Trelat}. For each model, geometrical studies have been proposed to analyze the controllability properties and the structure of extremals. For the sake of simplicity, we restrict our analysis to coplanar transfers at constant mass. The satellite dynamic is described by equations expressed in radial-orthoradial coordinates $(P,e_x,e_y,L)$ as follows:
\begin{eqnarray}
\left[\begin{array}{c}
\dot{P}\\
\dot{e}_x\\
\dot{e}_y\\
\dot{L}
\end{array}\right] = \sqrt{\frac{\mu}{P}}\left[\begin{array}{c}
0\\
0\\
0\\
W^2/P
\end{array}\right]
 + \sqrt{\frac{P}{\mu}}\left[\begin{array}{c}
0\\
sin~L\\
-cos~L\\
0
\end{array}\right]u_1\nonumber\\ +\sqrt{\frac{P}{\mu}}\left[\begin{array}{c}
2P/W\\
cos~L +(e_x+cos~L)/W\\
sin~L +(e_y+sin~L)/W\\
0
\end{array}\right]u_2\label{satellite:system}
\end{eqnarray}
where: $W=1+e_x cos~L+e_y sin~L$. Variables $(P,e)$ are respectively the parameter and the eccentricity vector of the osculating ellipse and $L$ the longitude. The variable change between the Cartesian and radial-orthoradial frames is defined by:
\begin{equation}
\left\{\begin{array}{lcl}
r_1 &=& \displaystyle\frac{P}{W}cos~L\\
\\
r_2 &=& \displaystyle\frac{P}{W}sin~L
\end{array}\right.\quad \quad\left\{\begin{array}{lcl}
v_1 &=& -\displaystyle\sqrt{\frac{\mu}{P}}(e_y+sin~L)\\
\\
v_2 &=& \displaystyle\sqrt{\frac{\mu}{P}}(e_x+cos~L)
\end{array}\right.\label{sat:chgt:repere}
\end{equation}
where $v=\dot{r}$ denotes the satellite velocity. From now on, the control constraints are expressed as:
\begin{equation*}
\abs{u_1}+\abs{u_2}<F_{max}.
\end{equation*}
and the control domain $\U_m$ could be defined by: $\U_m=Conv\left([-F_{max},0],[0,-F_{max}],[F_{max},0],[0,F_{max}] \right).$\\

The control problem is to minimize the time transfer of the satellite from a given initial orbit to a given final one, the position of the satellite on this last orbit being free. The related initial and final conditions are summed up on Table \ref{tab:cd}.
\begin{table}[htbp]
\begin{center}
\begin{tabular}{|ll|}
\hline
{\bf Initial conditions} & {\bf Final conditions}\\
\hline
$P^0=11.625$ Mm & $P^f=42.165$ Mm\\
$e_x^0=0.75$ & $e_x^f=0$\\
$e_y^0=0$ & $e_y^f=0$\\
$L^0=\pi$ rad & $L^f$ free\\
\hline
\end{tabular}
\caption{Initial and final condition for the orbital transfer problem}\label{tab:cd}
\end{center}
\end{table}

\begin{figure}[htbp]
\begin{center}
\scalebox{0.35}{\includegraphics{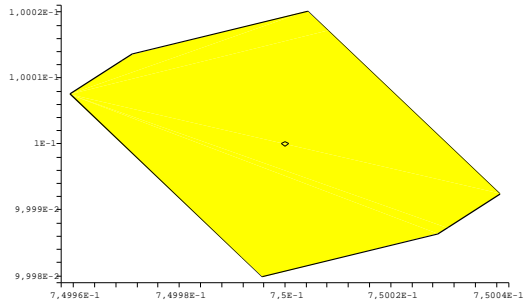}}\\
(a)\\
\scalebox{0.5}{\includegraphics{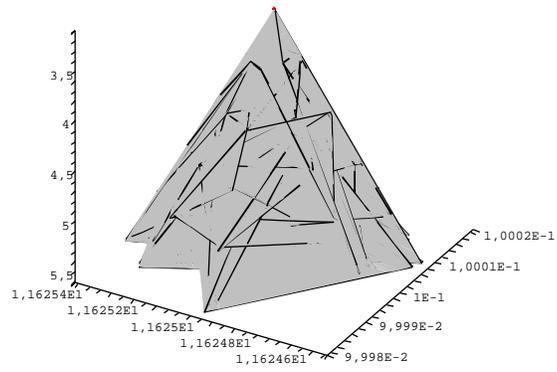}}\\
(b)
\caption{Attainable sets for (a) the eccentricity $e=(e_x,e_y)$ (b) parameters $(P,e_y,L)$, at the rougher scale ($h=6$).}\label{fig:atteignable:satellite:1}
\end{center}
\end{figure}

Let us now consider an initial point $X_0=[11.625, 0.75, 0, \pi]$ expressed in orbital coordinates. At the rougher scale $h=6$, our algorithm \ref{algo:convex:approximation} computes 451 controllable points in 4 minutes, building a four-dimensional approximation made of 104 vertices. By projection, we could obtain for example:
\begin{itemize}
\item the attainable set of the eccentricity $e$ from the initial point $e_0=(0.75,0)$ (see figure \ref{fig:atteignable:satellite:1}-(a)).
\item the three dimensional attainable set of parameters $(P,e_y,L)$ (see figure \ref{fig:atteignable:satellite:1}-(b)).
\end{itemize}
\begin{figure}[htbp]
\begin{center}
\scalebox{0.4}{\includegraphics{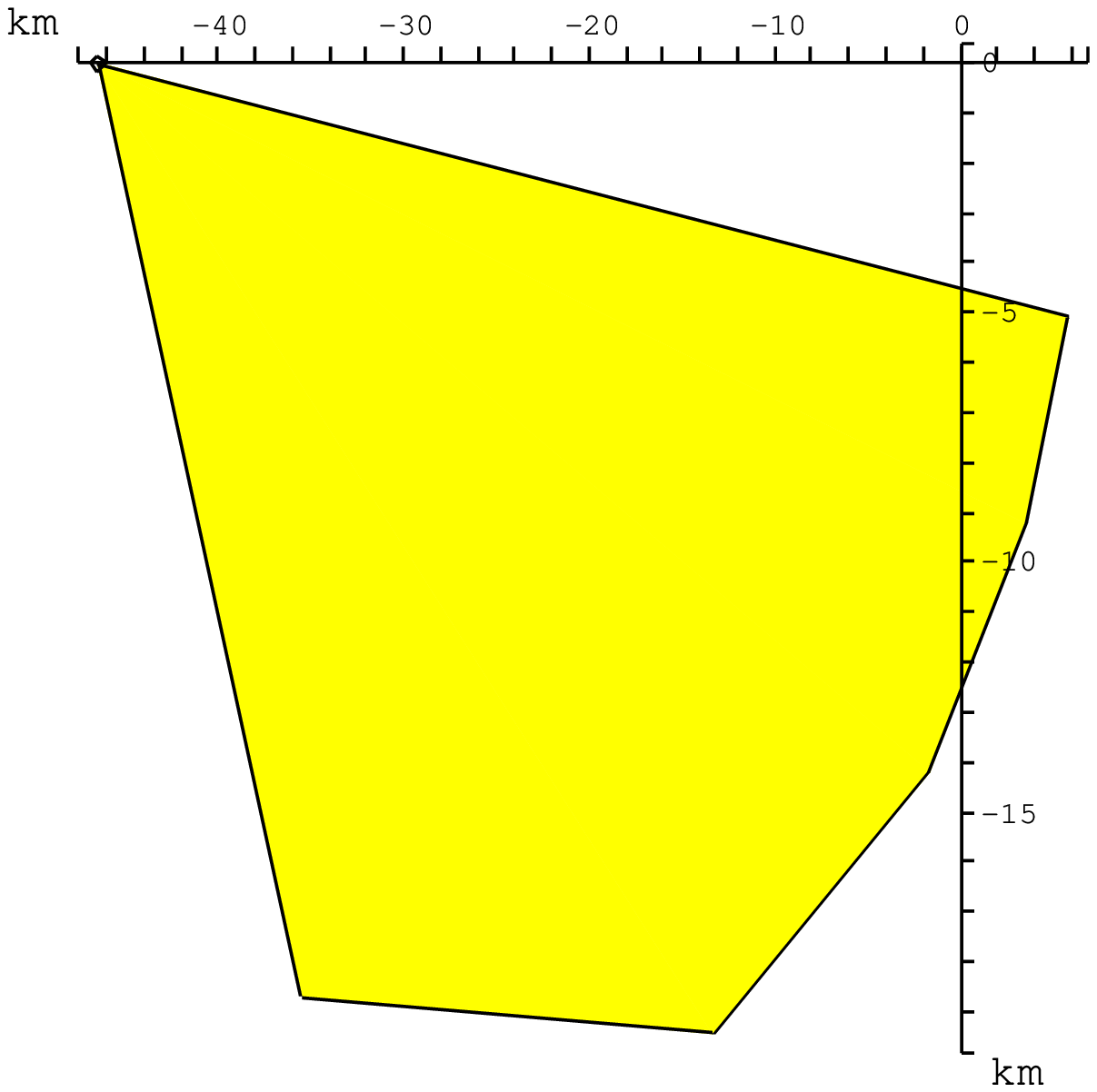}}
\caption{Approximation of the satellite attainable domain at scale $h=6$ in the Cartesian frame $(r_1,r_2)$}\label{fig:atteignable:satellite:2}
\end{center}
\end{figure}
Using variable change (\ref{sat:chgt:repere}), these numerical results could be translated in Cartesian coordinates $(r_1,r_2)$ and give an approximation of the attainable set from $X_0$ into the Earth Cartesian frame (see figure \ref{fig:atteignable:satellite:2}). 

Observe that since the Earth radius $r_{\cal E}$is about 6 kilometers, collision with the Earth is possible. To avoid this configuration, we could add a constraint like $\|r\|>r_{\cal E}$: we would then have to redefine the involved cells of the state mesh and this way, to restrict the controllable domain.

\section{Approaching nonlinear optimal solutions by optimal solutions of the hybrid control problem}\label{sec:optimality}
%
%
In this section we focus on the resolution of the nonlinear optimal control problem via the hybrid approach. Let $\U_m$ be a convex and compact polytope in $\R^m$ such that: $0\in \U_m$: in section \ref{sec:hybrid:automaton}, we shaw that any nonlinear control system like (\ref{syst:nonlinear}) can be approached by a piecewise affine hybrid model. This enables us to define two optimal control problems related to each control system:

\begin{prob}{NL}{
Minimize the cost function $J(X_0,u) = \int_0^{t_f} l(X(t),u(t))dt$ with respect to the control $u$ under the nonlinear dynamic:
$$\left\{
\begin{array}{l}
\dot{X}(t) = f(X(t),u(t))\\
X(0) = X_0
\end{array}
\right.~~~~~~~~~~
X(t_f) \in \target_f
$$
and the constraints: $\forall t \in [0,t_f], u(t) \in \U_m$. The final time $t_f$ is unspecified.}
\end{prob}

\begin{prob}{{\cal H}}{
Minimize the cost function $J(X_0,u) = \int_0^{t_f} l(X(t),u(t))dt$ with respect to the control $u$ under the piecewise affine hybrid dynamic:
$$\left\{
\begin{array}{l}
\dot{X}(t) = f_h(X(t),u(t))\\
X(0) = X_0
\end{array}
\right.~~~~~~~~~~
X(t_f) \in \target_f
$$
and the constraints: $\forall t \in [0,t_f], u(t) \in \U_m$. The final time $t_f$ is unspecified.}
\end{prob}

As for the controllability problem, a key point of our hybrid approach is first to prove that solutions of the hybrid problem \pb{\cal H} are satisfactory approximations of those of \pb{NL}. We then address the problem of solving the hybrid optimal control problem \pb{\cal H}: a powerful tool for optimal trajectories analysis is the (PMP) Pontryagin maximum principle \cite{Pontryagin}. First, in section \ref{ssec:hybrid:PMP}, we state a suitable maximum principle for our hybrid optimal control problem \pb{\cal H} and then deduce in section \ref{ssec:hyb:trajectories} the general structure of hybrid optimal trajectories.


\subsection{Hybrid maximum principle}\label{ssec:hybrid:PMP}
In the classical context of smooth systems, the Pontryagin maximum principle provides first-order necessary optimality conditions and a pseudo-Hamiltonian formulation of the considered optimal control problem. It has been first stated by Pontryagin and al. in \cite{Pontryagin}. Stronger versions have been stated especially in the context on constraints on the control. References are numerous~; one can refer to \cite{Pontryagin,Lee:Markus,Clarke,Agrachev:Sachkov} for a proof of the maximum principle or to \cite[paragraphe 1.5]{Bonnard:Caillau} and \cite{Trelat:book} for a detailed study of its applications.

The main difficulty of establishing a hybrid maximum principle is to manage discontinuities generated by the discrete component of the hybrid system. Today several formulations has been proposed depending on the structure and the regularity of the studied hybrid systems \cite{Sussmann:99b,Zaytoon}. Especially P. Riedinger, C. Iung and F. Kratz state and prove a hybrid maximum principle for a general class of hybrid systems in \cite[theorem 2.3]{Riedinger:Iung:Kratz:03}. Next our hybrid maximum principle is inspired from their theorem.

\begin{theorem}[Hybrid maximum principle]
We introduce the Hamiltonian $H_h$ related to the hybrid model ${\cal H}$ defined in section \ref{sec:hybrid:automaton}:
$$H_h(X,u,\lambda) = l(X,u)+\lambda^\top f_h(X,u).$$
If the control $u$ is optimal on the time interval $[0,t_f]$, then there exists a piecewise absolute continuous application $\lambda:[0,t_f]\rightarrow \R^n$ such that:
\begin{itemize}
\item[$\bullet$] For almost all $t\in [0,t_f]$, $H$ reaches its minimum with regards to the control at point $u(t)$ i.e.: $$H_h(X(t),u(t),\lambda(t)) = \min\limits_{v\in \U_m} H_h(X(t),v,\lambda(t)).$$
\item[$\bullet$] $\dot{X}(t) = \displaystyle\frac{\partial H_h}{\partial \lambda}(X(t),u(t),\lambda(t))$
\item[$\bullet$] $\dot{\lambda}(t)^\top = -\displaystyle\frac{\partial H_h}{\partial X}(X(t),u(t),\lambda(t))$  (Euler-Lagrange equation)
\item[$\bullet$] $H_h(X(t),u(t),\lambda(t))=0$ along the optimal trajectory.
\item[$\bullet$] At each transition instant $t_i$ between two modes $q_i$ and $q_{i+1}$, we have the following transversality conditions:
$$\left\{\begin{array}{l}
\lambda(t_i^+)=\lambda(t_i^-)\\
\\
H_h(X(t_i^+),u(t_i^+),\lambda(t_i^+)) = H_h(X(t_i^-),u(t_i^-),\lambda(t_i^-))
\end{array}\right.$$
\end{itemize}\label{th:hybrid:maximum}
\end{theorem}

We first recall that by definition, the state variable is continuous at each transition time $t_i$: $X(t_i^+)=X(t_i^-)$, which together with the autonomous feature of the hybrid system, implies the continuity of the hybrid Hamiltonian at each transition time.

Then transitions between discrete modes of our hybrid automaton are non constrained: at each transition time $t_i$ between (at least) two modes $q_i$ and $q_{i+1}$, our hybrid system could either switch to mode $q_{i+1}$, or stay in mode $q_i$, depending on the value of the field vector $f_h$ at $X(t_i)$. This implies the continuity of the adjoint vector $\lambda$ at each transition time.

\subsection{Characterization of hybrid optimal trajectories}\label{ssec:hyb:trajectories}

\subsubsection{Necessary optimality conditions in mode $q$}\label{ssec:nec:CD}
Let $X_0\in \R^n$ be a controllable point by our hybrid system ${\cal
  H}$ and $q$ the corresponding discrete mode satisfying: $X_0\in
D_q$. We assume that the initial time $t=0$ is not a transition time
of the system, which means that $X_0$ does not belong to any guard of
${\cal H}$ (in other words $X_0\in int~D_q$). This assumption guarantees that any (state) trajectory steered from $X_0$ remains inside the cell $D_q$ on a non empty time interval. 

Our point is now the local resolution (inside $D_q$) of the optimization problem stated by the hybrid maximum principle. According to this principle, any state and control trajectory candidate to optimality, also called extremal, in mode $q$ has to fulfill three conditions: to minimize the hybrid Hamiltonian, to enable the definition of an appropriate adjoint state vector $\lambda$ and to cancel the Hamiltonian along this trajectory.

\paragraph{Minimization of the hybrid Hamiltonian $H_h$ in mode $q$} 
In paragraph \ref{ssec:local:data}, we defined the set of local control constraints at a given position $X$ in one mode $q$ by introducing the sets:
$$\forall X\in D_q,\forall q'\in {\cal K}(q),U_{q'}(X)=\{u\in\U_m;(X,u)\in\Delta_{q'}\}.$$
where ${\cal K}(q)$ denotes the set of indices of state-control cells of the mesh $\Delta$ whose projection is exactly $D_q$. Therefore $(U_{q'}(X))_{q'\in {\cal K}(q)}$ forms a partition of the control domain $\U_m$ i.e.: 
\begin{equation}
\forall X\in D_q, \U_m = \bigcup\limits_{q'\in{\cal K}(q)} U_{q'}(X)\label{partition:Um}
\end{equation}
Now, by definition of the hybrid system $\cal H$, between two successive transition instants $t_i$ and $t_{i+1}$, the state trajectory $X(t)$ remains in cell $D_q$, the state and control trajectory $(X(t),v)$ belongs to $\Delta_{q'}$ for all $v\in U_{q'}(X(t))$ and the hybrid dynamic at any point $(X(t),v)$ for $v\in U_{q'}(X(t))$ is affine: $\forall t\in [t_i,t_{i+1}],\forall v\in U_{q'}(X(t)),f_h(X(t),u(t)) = A_{q'}X(t)+B_{q'}v+c_{q'}$. Together with (\ref{partition:Um}), we then have:
$$\begin{array}{l}
\inf\limits_{v\in \U_m} H(X(t),v,\lambda(t)) \\
\phantom{blabla}= \min\limits_{q'\in {\cal K}(q)} \left(\min\limits_{v\in U_{q'}(X(t))} H(X(t),v,\lambda(t))\right)\\
\phantom{blabla}= \min\limits_{q'\in {\cal K}(q)} \left(\min\limits_{v\in U_{q'}(X(t))} l(X(t))+\lambda(t)^\top f_h(X(t),v)\right)\\
\phantom{blabla}= \min\limits_{q'\in {\cal K}(q)} \left(\min\limits_{v\in U_{q'}(X(t))} l(X(t))+\lambda(t)^\top A_{q'}X(t)\right.\\
\left.\phantom{\displaystyle\min\limits_{v\in U_{q'}(X(t))}blablablabablablabla}+\lambda(t)^\top c_{q'}+\lambda(t)^\top B_{q'}v\right)
\end{array}$$
Let us consider the inner minimization problem for a given $q'\in {\cal K}(q)$. The objective is affine in the optimization variable $v$. Moreover, according to the proposition \ref{prop:Uq}, the set $U_{q'}(X(t))$ is a $m$-simplex whose vertices are denoted by $\sigma_i(t)=F_iX(t)+g_i$, $i=1,\dots,m+1$. The induced constraints on the control $v$ are therefore affine and the inner minimization problem actually is a linear program whose solution is well-known:
$$\begin{array}{c}
v^\star(t) = F_iX(t)+g_i \\
\mbox{ if }~\forall j\neq i, \lambda(t)^\top B_{q'}\left[(F_i-F_j)X(t)+g_i-g_j\right]<0.
\end{array}$$
We observe that the switching conditions explicitly depend both on the state $X$ and the adjoint state $\lambda$.

\begin{remark}[Switching curves and singular controls]
From now on, we assume that optimal controls are non singular or equivalently, that zeros of each switching function 
$$S_{i,j}:t\mapsto \lambda(t)^\top B_{q'}\left[(F_i-F_j)X(t)+g_i-g_j\right]$$
where $i,j=1,\dots,m,i\neq j$, are isolated.
\end{remark}

Consequently our minimization problem in mode $q$ could be replaced by
the simultaneous solving of a finite number (exactly $card~{\cal
  K}(q)$) of linear programs: the searched optimal controls are
feedback (more precisely they are piecewise affine in the state variable $X$). Notice that the choice of each piece $FX+g$ also depends on the adjoint state variable $\lambda$. \\

The next step is to check that such controls fulfill the missing optimality necessary conditions given by the hybrid maximum principle. Let us now introduce a given affine feedback control $u=FX+g$ and the resulting trajectory in mode $q$.
\begin{figure}[htbp]
\begin{center}
\begin{minipage}{.49\textwidth}
\begin{center}
\scalebox{.5}{\input{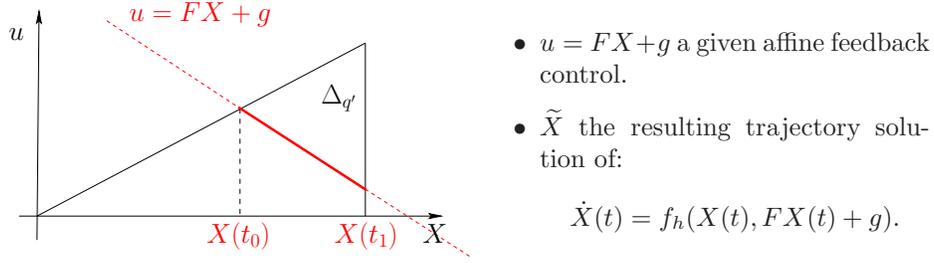}} 
\end{center}
\end{minipage}
\begin{minipage}{.5\textwidth}
\begin{itemize}
\item $u=FX+g$ a given affine feedback control.
\item $\widetilde X$ the resulting trajectory solution of: 
$$\dot{X}(t) = f_h(X(t),FX(t)+g).$$
\end{itemize}
\end{minipage}
\caption{Definition of an affine feedback control $u=FX+g$ in mode $q$}\label{fig:CNoptim}
\end{center}
\end{figure}
As illustrated on figure \ref{fig:CNoptim}, we now focus on a time interval $[t_0,t_1]$ on which the considered trajectory $(\widetilde X,F\widetilde X+g)$ remains in a state-control cell $\Delta_{q'}$ for some $q'\in {\cal K}(q)$:
$$\forall t\in [t_0,t_1],(\widetilde X(t),F\widetilde X(t)+g)\in\Delta_{q'}.$$
It follows that the hybrid dynamic is affine along the trajectory $(\widetilde X,F\widetilde X+g)\in \Delta_{q'}$ i.e. for all $t\in [t_0,t_1]$: 
$$f_h(\widetilde X(t),F\widetilde X(t)+g) = \left(A_{q'}+B_{q'}F\right)\widetilde X(t)+B_{q'}g+c_{q'}.$$

\paragraph{Definition of the adjoint vector $\lambda$ in mode $q$}
The adjoint parameter is then simply defined as solution of the following first order differential equation:
\begin{equation}
\dot{\lambda}(t)^\top = -\frac{\partial l}{\partial X}(\widetilde X(t)) - \lambda(t)^\top(A+BF)\label{def:lambda}
\end{equation}
for all $t\in ]t_0,t_1[$. According to our regularity assumptions, the cost function $l$ is continuously differentiable on $\R^n\times \U_m$, so that $t\mapsto \frac{\partial l}{\partial X}(\widetilde X(t))$ is continuous. This implies the expected absolute continuity of the solution $\lambda$. Notice that no initial condition on $\lambda$ is required.

\paragraph{Hybrid Hamiltonian equal to $0$ along the optimal trajectory}
Lastly we have to check that the hybrid Hamiltonian is equal to $0$ along the considered trajectory $(\widetilde X,F\widetilde X+g)$. We write for all $t\in [t_0,t_1]$:
$$\begin{array}{l}
H(\widetilde X(t),F\widetilde X(t)+g,\lambda(t))\\
\phantom{blablabla}= l(\widetilde X(t))+\lambda(t)^\top\left[(A+BF)\widetilde X(t)+Bg+c\right]\\
\phantom{blablabla}= l(\widetilde X(t))+\lambda(t)^\top\dot{\widetilde X}(t).
\end{array}$$
Using $\ddot{X}(t)=(A+BF)\dot{X}(t)$ and relation (\ref{def:lambda}), we prove
$$\begin{array}{lcl}
\displaystyle\frac{d}{dt}H(X(t),u(t),\lambda(t)) \\
\phantom{blabla}= \displaystyle\frac{\partial l}{\partial X}(X(t))\dot{X}(t) + \dot{\lambda}(t)^\top\dot{X}(t) +\lambda(t)^\top\ddot{X}(t)\vspace{.2cm}\\
\phantom{blabla}= \displaystyle\frac{\partial l}{\partial X}(X(t))\dot{X}(t)+\lambda(t)^\top(A+BF)\dot{X}(t)\\
\phantom{blabla}\phantom{blabla} -\left(\displaystyle\frac{\partial l}{\partial X}(X(t))+\lambda(t)^\top(A+BF)\right)\dot{X}(t) \\
\phantom{blabla}= 0
\end{array}$$
for all $t\in [t_0,t_1]$. That means that the Hamiltonian is constant along the considered trajectory $(\widetilde X,F\widetilde X+g,\lambda)$. Thus to fulfill the last optimality necessary condition required by our hybrid maximum principle, we just have to cancel the hybrid Hamiltonian at the initial time i.e.:
\begin{equation}
l(X(t_0))+\lambda(t_0)^\top\left[(A+BF)X(t_0) + Bg+c\right]=0.\label{cd:local:optimality}
\end{equation}
This gives an explicit constraint on the adjoint state variable $\lambda$ at the initial time $t_0$.

In conclusion that feedback affine controls are good local candidates to optimality. Next we propose a strategy to globally reconstruct the extremals of our hybrid optimal control problem \pb{{\cal H}}.

\subsubsection{Global structure of hybrid extremals}\label{ssec:extremal}
Let $(X,u,\lambda)$ be an extremal of the hybrid optimal control problem \pb{{\cal H}} for the initial condition $X(0)=X_0$. By definition \ref{def:hybrid:execution} of an hybrid trajectory, there exists a sequence of time intervals $([t_i,t_{i+1}])_{i=0\dots r}$ and a sequence of discrete modes $(q_i)_{i=0\dots r+1}$ such that:
\begin{equation}
\forall t\in [t_i,t_{i+1}], X(t) \in D_{q_i}
\end{equation}
The transversality conditions stated by the hybrid maximum principle
together with the definition of our hybrid model, guarantee the
continuity of both the state $X$ and the adjoint state $\lambda$ at
each switching time $t_i$ between two discrete modes. There just remains to determine the structure of trajectories inside each mode $q_i$.

\begin{figure}[htbp]
\begin{center}
\scalebox{.55}{\input{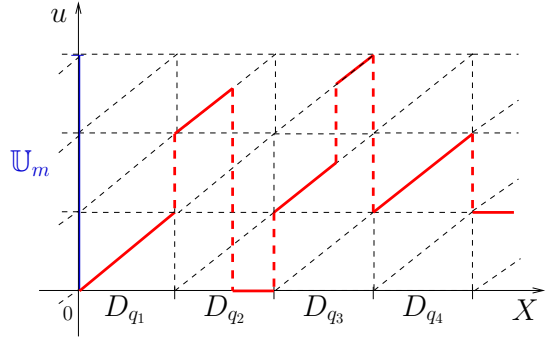}} 
\caption{Hybrid extremals evolve along edges of the state and control mesh $\Delta$}\label{fig:trajectoire:arete}
\end{center}
\end{figure}

\paragraph{Structure in mode $q_i$} We focus on one discrete mode $q_i$ of our hybrid automaton ${\cal H}$ and introduce the set $\Sigma_{q_i}(X)$ of the vertices of the control domain $U_{q'}(X)$ induced at position $X$ in mode $q_i$:
\begin{equation*}
\begin{array}{lcl}
\Sigma_{q_i}(X) &=& \{FX+g \in \U_m / \exists q'\in {\cal K}(q_i),\\ 
&&(FX+g) \mbox{ is a vertex of } U_{q'}(X)\}\subset \U_m.
\end{array}\label{sigma}
\end{equation*}
According to results presented in paragraph \ref{ssec:nec:CD}, the
optimal control $u$ is locally the solution of a linear program and can be expressed as an explicit piecewise affine function of the state $X$, each piece's choice implicitly depending on the adjoint vector $\lambda$: therefore there exists a subdivision $(\widetilde t_j)_{j=0\dots N_i}$ of the time interval $]t_i,t_{i+1}[$ such that for all $j\in\{0,\dots,N_i-1\}$:
\begin{equation}
\begin{array}{c}
\exists (F_j,g_j)\in \M_{m,n}\times \R^m, \forall t\in [\widetilde t_j,\widetilde t_{j+1}[,\\
u(t) = F_jX(t)+g_j \in \Sigma_{q_i}(X(t)).
\end{array}\label{local:optimal:control}
\end{equation}
Hence: $(X,F_jX+g_j)\in \Delta_{q'}$ on $[\widetilde t_j,\widetilde t_{j+1}[$. In other words, as illustrated on figure \ref{fig:trajectoire:arete}, the optimal trajectory $(X,u)$ evolves along the edges of the state and control mesh $\Delta$, continuously in $X$ and not in $u$.

The main issue is so to efficiently compute both the optimal control
and the adjoint state variable $\lambda$ at any given position
$X$. Indeed Euler-Lagrange equations determining $\lambda$ depend on
the control choice, which in turns depends on $\lambda$. Therefore we
need to determine $u=FX+g$ and $\lambda$ satisfying the optimality
conditions without having to compute any trajectories.

\paragraph{Switching conditions and optimal control computation}
Let us consider our hybrid optimal control problem at a given time
$\widetilde t_j$. The state $\widetilde X_j=X(\widetilde t_j)$ and the
adjoint state $\widetilde \lambda_j=\lambda(\widetilde t_j)$ are
supposed to be known. According to (\ref{local:optimal:control}), we
want to determine control candidate to optimality satisfying $u=FX+g$
in the set $\Sigma_{q_i}(X)$. We also need to compute the first time $\widetilde t_{j+1}>\widetilde t_j$ at which the optimality conditions are not fulfilled any more.

We assume that the control $u=FX+g$ satisfies all the optimality conditions required by the hybrid maximum principle at time $\widetilde t_j$: 
$$\begin{array}{lcl}
H_h(\widetilde X_j,F\widetilde X_j+g,\widetilde\lambda_j)&=&\min\limits_{u\in \U_m} H_h(\widetilde X_j,u,\widetilde\lambda_j) \nonumber\\
&=& \displaystyle\min_{\scriptsize\begin{array}{c}
F'\widetilde X_j+g'\\
\in \Sigma_{q_i}(\widetilde X_j)
\end{array}} H_h(\widetilde X_j,F'\widetilde X_j+g',\widetilde \lambda_j)\nonumber\\
H_h(\widetilde X_j,F\widetilde X_j+g,\widetilde\lambda_j)&=&0\nonumber
\end{array}$$
This induces that switching functions could be expressed via the hybrid Hamiltonian this way:
$$\begin{array}{l}
\displaystyle S_{FX+g,F'X+g'} = \lambda^\top\left[f_h(X,FX+g)-f_h(X,F'X+g')\right]\\
\phantom{\displaystyle S_{FX+g,F'X+g'}} = -l(X) - \lambda^\top f_h(X,F'X+g')\\
\phantom{\displaystyle S_{FX+g,F'X+g'}} = -H_h(X,F'X+g',\lambda)\leq 0
\end{array}$$
Therefore switching conditions are obtained by studying the zeros of
the application: $u\mapsto H_h(\widetilde X_j,u,\widetilde\lambda_j)$
over the finite set $\Sigma_{q_i}(\widetilde X_j)$. We thus compute:
$H_h(\widetilde X_j,F\widetilde X_j+g,\widetilde\lambda_j)$ for all
$F\widetilde X_j+g \in \Sigma_{q_i}(\widetilde X_j)$. When the latter equals zero, the corresponding control is a good candidate to optimality. The main advantage of this approach is that we only need $card~\Sigma_{q_i}(\widetilde X_j)$, i.e. at most $m+card~{\cal K}(q_i)$, function evaluations of $H_h(\widetilde X_j,\cdot,\widetilde\lambda_j)$.

We then introduce:
$$\begin{array}{l}
\widetilde t_{j+1} = \min\left(t_{i+1};\right.\\
\phantom{\widetilde t_{j+1} = \min(}\left.\inf\{ t>t_j;H_h(\widetilde X(t),F\widetilde X(t)+g,\widetilde \lambda(t))>0\}\right)\end{array}$$
where ($\widetilde X,\widetilde \lambda)$ denotes the state and adjoint state trajectory steered from ($\widetilde X_j,\widetilde \lambda_j)$ at time $\widetilde t_j$ according to the control $u=FX+g$. In other words, $\widetilde t_{j+1}$ denotes either the exit time, namely $t_{i+1}$, from the current mode $q_i$ or the first instant at which optimality conditions are not satisfied anymore. At time $\widetilde t_{j+1}$, we have new initial conditions $(\widetilde X_{j+1},\widetilde \lambda_{j+1})=(\widetilde X(\widetilde t_{j+1}),\widetilde \lambda(\widetilde t_{j+1}))$.

\subsection{Simulation algorithm of hybrid extremals}
For a given initial condition of the adjoint state, we propose an
algorithm to compute extremals of the hybrid optimal control problem
i.e. trajectories satisfying the necessary optimality conditions of
the hybrid maximum principle. Algorithm \ref{algo:extremal} details
this computation of a hybrid extremal in a given mode $q_i$.
\begin{algorithm}[htbp]
\caption{\texttt{\bf Extremal computation in a given discrete mode $q_i$}}\label{algo:extremal}
\begin{algorithmic}[1]
\REQUIRE $t_i\geq 0$ initial time in mode $q_i$ at the position $X_i\in D_{q_i}$ for the given initial condition: $\lambda(t_i)=\lambda_i$.
\STATE{\bf Initialization.} $j=0$; $\widetilde t_0 = t_i$; $\widetilde X_0 = X_i$; $\widetilde \lambda_0=\lambda_i$;
\STATE{\bf Control computation.} Find the control $F\widetilde X_j+g\in \Sigma_{q_i}(\widetilde X_j)$ satisfying the necessary optimality conditions at time $\widetilde t_j$ i.e.:
\begin{eqnarray}
H_h(\widetilde X_j,F\widetilde X_j+g,\widetilde \lambda_j) &=& 0\nonumber\\
H_h(\widetilde X_j,F'\widetilde X_j+g',\widetilde \lambda_j)&>&0\nonumber
\end{eqnarray}
for all $F'\widetilde X_j+g' \in \Sigma_{q_i}(\widetilde X_j)\backslash\{F\widetilde X_j+g\}$.
\STATE{\bf State and adjoint state trajectories.} Compute the state $X_j = X[\widetilde X_j,FX+g]$ and adjoint state $\lambda_j=\lambda[\widetilde \lambda_j,FX+g]$ trajectories respectively steered from $\widetilde X_j$ and $\widetilde \lambda_j$ according to the control $u(X)=FX+g$:
$$\begin{array}{l}
\dot{X}(t) = f_h(X(t),FX(t)+g),~~X(t_i)=\widetilde X_j\\
\\
\begin{array}{r}
\dot{\lambda}(t)^\top = -\displaystyle\frac{\partial H_h}{\partial X}(X[\widetilde X_j,FX+g](t),\phantom{aaaaaaaaaa}\\
FX[\widetilde X_j,FX+g](t)+g,\lambda(t))\end{array}\\
\lambda(t_i)=\widetilde \lambda_j
\nonumber
\end{array}$$
\STATE Compute the exit time $t_{i+1}$ of $X[\widetilde X_j,FX+g]$ from the cell $D_{q_i}$ and the next cell $D_{q_{i+1}}$ at time $t_{i+1}$.
\STATE{\bf Switching time.} Introducing the Hamiltonian along the considered trajectory, compute 
$$\begin{array}{l}
\widetilde t_{j+1} = \inf\{t>t_j;H_h(X[\widetilde X_j,FX+g](t),\\
\phantom{\widetilde t_{j+1} =}FX[\widetilde X_j,FX+g](t)+g,\lambda[\widetilde \lambda_j,FX+g](t))>0\}\end{array}$$
\IF {$\widetilde t_{j+1}$ exists and $\widetilde t_{j+1}<t_{i+1}$}\STATE Go back to step 2 with new initial conditions: $\widetilde t_{j+1}$, $\widetilde X_{j+1}=X[\widetilde X_j,FX+g](\widetilde t_{j+1})$ and $\widetilde \lambda_{j+1}=\lambda[\widetilde \lambda_j,FX+g](\widetilde t_{j+1})$ until $\widetilde t_{j+1} = t_{i+1}$ i.e. until the state leaves the cell $D_{q_i}$.
\STATE $j=j+1$;
\ENDIF
\end{algorithmic}
\end{algorithm}

The global simulation algorithm of hybrid extremals relies on the computation of discrete modes reached by the trajectory satisfying the necessary optimality conditions and on the transversality conditions that guarantee the continuity of the adjoint state variable $\lambda$. The principle is as follows: starting from a given initial condition $(t_0,X_0,\lambda_0)$ in mode $q_0$ and using algorithm \ref{algo:extremal}, we compute the trajectory inside $D_{q_0}$ satisfying the hybrid maximum principle and its exit time from the cell $D_{q_0}$ towards an adjacent cell $D_{q_1}$. By continuity of the state and adjoint state, the exit conditions give new initial conditions $(t_1,X_1,\lambda_1)$ in mode $q_1$ and we pursue the simulation the same way.





\section{Conclusion}
In this paper, we address the problem of hybrid optimal control of
nonlinear dynamical systems. First we propose a hybrid approximation,
by way of a hybrid automaton with piecewise affine dynamics, of
complex systems. We therefore define a piecewise affine optimal control
problem and then focus on two main problems in control theory: the
controllability question and the computation of optimal solutions. 

For the controllability, we have given a convex approximation of the
controllable set which is at a guaranteed distance of the actual
controllable set. Furthermore, when the controllable set is convex,
this approximation is an under-approximation and thus guarantees
controllability. We then proposed and implemented an algorithm 
to compute this approximation and its efficiency and
accuracy are shown on the examples of the nonlinear spring and of the
change of orbit.

Then we addressed the question of optimality. By developing an hybrid
maximum principle we are able to characterize the hybrid optimal
trajectories and to prove that the optimal control have the form of
affine feedbacks. Moreover, the principle gives us a way to compute
the switching times between modes and thus an algorithm to follow the
optimal trajectories.

Further improvements include mixing this approach with shooting
methods in order to get a better initialization of the adjoint
vector.

\bibliographystyle{plain}
\bibliography{mybib}

\appendix

\section{Proof that the control constraints are $m$-simplices}\label{apdx:uq}
\begin{figure}[htbp]
\begin{center}
\includegraphics[width=.5\textwidth]{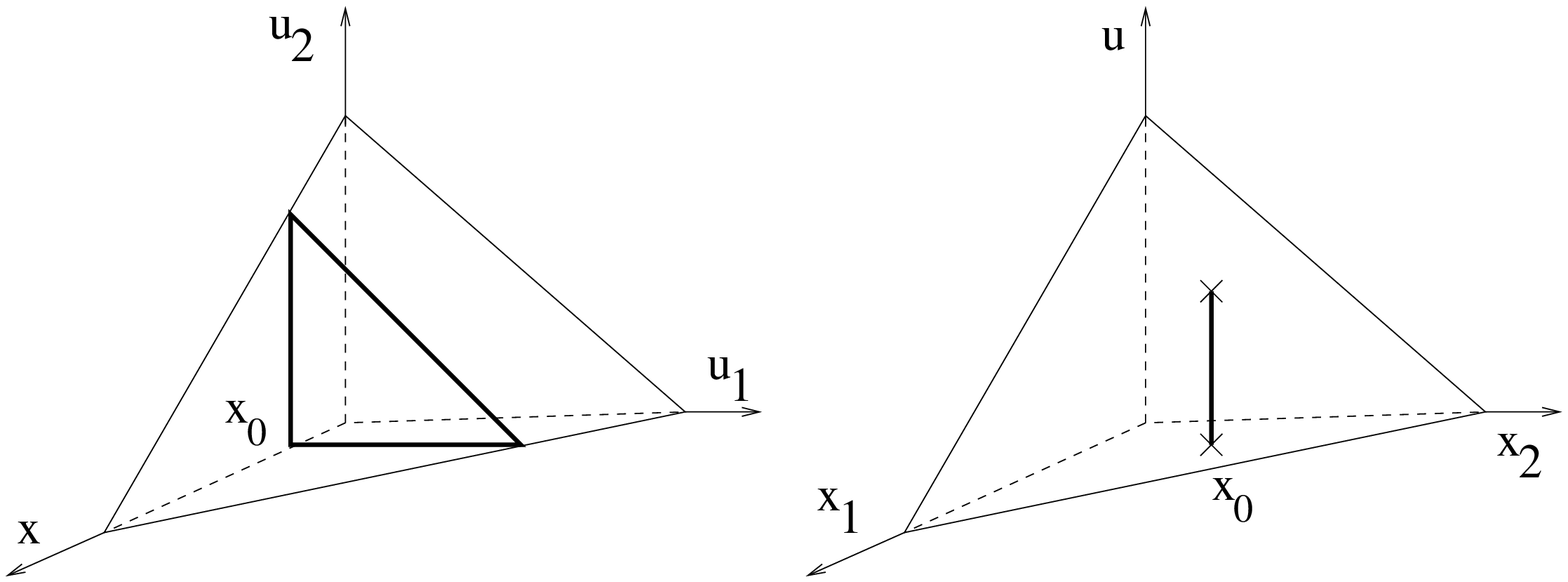}
\caption{Intersection $P\cap \Delta_{q'}$ for $(n,m)=(1,2)$ on the
left, $(n,m)=(2,1)$ on the right.}\label{fig:exemple:intersection}
\end{center}
\end{figure}

We take $q'\in{\cal K}(q)$. The problem is to compute the
intersection of the affine subspace $P=\{(X_0,u)/u\in \R^m\}$ in
$\R^{n+m}$ and the simplex $\Delta_{q'}$ (see figure
\ref{fig:exemple:intersection}):
$$U_{q'}(X_0) = \{(X_0,u)~/~u\in \R^m\} \cap \Delta_{q'} =P \cap \Delta_{q'}$$
We know that the intersection of an affine subspace with a simplex is a
simplex \cite[Theorem 1.1]{Ziegler}. Therefore $U_{q'}(X_0)$ is a simplex in
$\R^m$. There
remain to to compute its dimension. The idea is to
consider the intersection of the subspace $P$ of dimension $m$ with
the n-faces of $\Delta_{q'}$. These intersections, when they exist,
are of dimension $0$, and define the vertices of $U_{q'}(X_0)$. We
then have to prove that there exists exactly $(m+1)$ different n-faces
intersecting $P$ to conclude.\\

Let $\{X_1,\dots,X_{n+1}\}$ be the vertices of $D_q$. By construction,
we have: $p_{\R^n}^{\bot}(\Delta_{q'})=D_q$. We can then introduce the
vertices of $\Delta_{q'}$:
$$\begin{array}{l}
V_k=\left[\begin{array}{c}
X_{i_k}\\
u_{j_k}
\end{array}\right],i_k\in \{1,\dots,n+1\},j_k\in
\{1,\dots,m+1\},\\
k=1,\dots,n+m+1
\end{array}$$
and: $\{X_{i_k};k=1,\dots,n+m+1\}=\{X_1,\dots,X_{n+1}\}$.\\
We also assume: $X_0\in \mathring{D_q}$, hence:
\begin{equation}
\exists !(\alpha_k)_{k=1\dots n+1} \in
]0,1[^{n+1},\sum\limits_{k=1}^{n+1}\alpha_k=1 \mbox{ and } X_0=\sum\limits_{k=1}^{n+1}\alpha_k
X_k\label{rel:convexe}
\end{equation}

\paragraph{{\bf $1^{st}$ step: any n-face of $\Delta_{q'}$ has either
one intersection with $P$, or none.}}~\\
Let $F$ be a $n$-face of $\Delta_{q'}$ defined by the list (renumbered
to ease the readability) of its vertices: 
$(\left[\begin{array}{c}
X_{i_k}\\
u_{j_k}
\end{array}\right])_{k=1,\dots,n+1}$
where $i_k\in \{1,\dots,n+1\}$, $j_k\in \{1,\dots,m+1\}$.\\
\begin{itemize}
\item {\it $1^{st}$ case: the $i_k$ are all different.}\\ 
We then have: $\{X_{i_k}~/~k=1,\dots,n+1\}=\{X_1,\dots,X_{n+1}\}$.
(\ref{rel:convexe}) thus becomes:
$X_0=\sum\limits_{k=1}^{n+1}\alpha_{i_k} X_{i_k}$. We state: $X_F=\sum\limits_{k=1}^{n+1} \alpha_{i_k}\left[\begin{array}{c}
X_{i_k}\\
u_{j_k}
\end{array}\right]=\left[\begin{array}{c}
X_0\\
\sum\limits_{k=1}^{n+1} \alpha_{i_k} u_{j_k}
\end{array}\right]$.

Then: $X_F\in P$. Moreover, by construction of $X_F$, $X_0\in \mathring{D_q}$ involves: $X_F \in \mathring{F}$, hence:  $X_F\in P\cap F$.\\
By uniqueness of the convex decomposition in a given simplex, we then
deduce:
$$X_F \in \mathring{F} \mbox{ and } P\cap F=\{X_F\} $$
\item {\it $2^{nd}$ case: the $i_k$ are not all different.}~\\
We then deduce: $\{X_{i_k}~/~k=1,\dots,n+1\} \subsetneq
\{X_1,\dots,X_{n+1}\}$. For $Y\in F$, we therefore have:
$$\exists ! (\gamma_{i_k})_{k=1\dots n+1} \in
[0,1]^{n+1},
\left\{\begin{array}{l}
\sum\limits_{k=1}^{n+1}\gamma_{i_k}=1\\
Y=\sum\limits_{k=1}^{n+1}\gamma_{i_k} \left[\begin{array}{c}
X_{i_k}\\
u_{j_k}
\end{array}\right]
\end{array}\right.$$
However, according to the hypothesis: $X_0\in \mathring{D_q}$, $X_0$ depends
on all the $X_i$, $i=1,\dots,n+1$. It follows:
$\sum\limits_{k=1}^{n+1}\gamma_{i_k}X_{i_k} \neq X_0$, i.e.: $Y\notin
P$, and: $F\cap P=\emptyset$.
\end{itemize}
\paragraph{{\bf $2^{nd}$ step: there exist exactly $(m+1)$ n-faces
intersecting $P$.}}~\\
By construction, we know that:
$\{X_{i_k}~/~k=1,\dots,n+m+1\}=\{X_1,\dots,X_{n+1}\}$ ; after
renumbering the vertices of $\Delta_{q'}$, we can assume:
$$\forall l \in \{1,\dots,n+1\},~V_l=\left[\begin{array}{c}
X_l\\
u_{j_l}
\end{array}\right]$$
\begin{itemize}
\item[$\bullet$] We state: $F_0=Conv(\{V_1,\dots,V_{n+1}\})$. By
affine independence of $X_1,\dots,X_{n+1}$ of $D_q$, the
vertices $V_1,\dots,V_{n+1}$ are also affinely independent, so that
$F_0$ is actually a n-face of $\Delta_{q'}$ intersecting $P$.

\item[$\bullet$] For $n+2\leq l \leq n+m+1$ (i.e. $m$ possible values
for $l$), we state:
$$\begin{array}{l}
F_l=Conv\left((\{V_1,\dots,V_{n+1}\}-\{V_{i_l}\})\cup \{V_l\}\right),\\
1\leq i_l\leq n+1
\end{array}$$
where: $V_l=(X_{i_l},u_{j_l})$ and $V_{i_l}=(X_{i_l},u_{j_{i_l}})$. As
for $F_0$, $F_l$ is a n-face of $\Delta_{q'}$ which intersects $P$.
\end{itemize}
To conclude, we have found $(1+m)$ n-faces of $\Delta_{q'}$ which
intersect $P$, i.e. $(m+1)$ vertices of $U_{q'}(X_0)$.

\end{document}